# A Relaxation-Based Decomposition Approach for Solving a Supported-Evacuation Problem in Wildfires

Shahryar Moradi[a,*], Antoine Sauré[a], Jonathan Patrick[a]

[a]*Telfer School of Management, University of Ottawa, 55 Laurier Ave., Ottawa, ON, K1N 6N5, Canada*

**Abstract**

Wildfires are among the most destructive natural disasters globally, particularly in North America, causing severe harm to communities, wildlife, and the environment. Given their growing threat, advanced disaster management strategies—especially concerning evacuation planning—have become necessary. Although extensive research exists on self-evacuation, limited attention has been paid to supported-evacuation, which involves assisting individuals unable to evacuate independently such as hospital patients, long-term care residents, and people with disabilities. To address this gap, we adopt a two-stage stochastic optimization approach that considers various types of patients, vehicles and medical facilities as well as time windows. We aim to optimize facility location, fleet sizing, and vehicle routing decisions to evacuate as many affected individuals as possible within time windows in a cost-efficient manner. We formulate two optimization models differing on whether the evacuation vehicles can be managed independently due to potential partial route disruptions. Given the NP-hard nature of the problem, we propose an innovative solution methodology based on Logic-Based Benders Decomposition, which also involves Combinatorial Benders Cuts forming a neighbourhood search and logic-based inequalities. The proposed methodology generates highly-quality solutions for realistically-sized instances within reasonable time frames. To highlight the benefits of using the proposed methodology, we compare its performance with that of three alternative evacuation policies through extensive numerical experiments. The results demonstrate significant improvements in shelter location decisions, vehicle usage, and evacuation times. Data collected during a community wildfire drill performed in Roxborough Park, Colorado (USA) are used to evaluate the performance provided by the two models in realistic conditions.

*Keywords:* OR in Disaster Relief, Supported-Evacuation in Wildfires, Facility Location and Vehicle Routing, Stochastic Programming, Problem Decomposition

## 1. Introduction

Wildfires, especially those close to residential areas, pose significant threats to both humans and wildlife. Despite extensive efforts to prevent them, catastrophic wildfires continue to occur all around the world (Zhou & Erdogan 2019) and represent one of the most common natural disasters in countries such as Canada, Chile, USA, and Australia. According to the Canadian National Fire Database, a yearly average of 8,000 wildfires occurred in Canada between 1970 and 2017, burning approximately 2.1 million hectares annually (National Forestry Database 2024). The 2016 Horse River Wildfire in northeastern Alberta is considered the costliest natural disaster in Canadian history, with damages amounting to $3.84 billion dollars (Tymstra et al.

*Corresponding author

*Email addresses:* Moradi@telfer.uottawa.ca (Shahryar Moradi), asaure@uottawa.ca (Antoine Sauré), Patrick@telfer.uottawa.ca (Jonathan Patrick)

2020). In Australia, wildfires caused on average \$2.5 billion dollars in damages annually between 1967 and 1999. During the same period, 223 wildfire-related deaths were recorded, representing 39% of all natural disaster-related fatalities in the country (Shahparvari & Abbasi 2017).

One of the most important challenges during wildfires is the timely evacuation of the affected population. There are numerous papers covering topics such as the timing of evacuation orders (McCaffrey et al. 2018, McLennan et al. 2019), evacuation behaviour (Toledo et al. 2018, Stasiewicz & Paveglio 2021), and wildfire and traffic simulation to assess the performance of alternative evacuation strategies (Beloglazov et al. 2016, Gwynne et al. 2023). However, these studies primarily focus on self-evacuation. The literature concerning the evacuation of those who cannot evacuate by themselves is rather scarce, despite the obvious significance it might have for people who need medical attention such as patients in hospitals and clinics, individuals with disabilities and limited mobility, seniors residing in long-term care facilities and retirement homes, people with no vehicles, and tourists visiting fire-effected areas (Kaiser et al. 2012).

Specifically, this paper concerns the evacuation of vulnerable people from fire-affected areas, where certain locations within the fire zone serve as temporary assembly areas housing or gathering patients who are unable to evacuate by themselves. To address this problem, we propose a novel two-stage stochastic (TSS) optimization approach that facilitates the design of timely and cost-efficient evacuation plans in the presence of several sources of uncertainty. The proposed approach allows for split pickups and split deliveries and optimizes the location of special-needs medical shelters (SMSs), vehicle requirements, and routing plans to evacuate as many individuals as possible in a timely and cost-efficient manner. Depending on whether or not vehicles of the same type can be managed independently, largely due to the possibility of partial route disruptions, we develop two optimization models. We determine conditions where even greedy approaches such as the heuristics that we use as benchmarks would have acceptable performance and discuss the possible benefits of cooperation between vehicles.

To be able to solve large-scale instances of the problem, we propose an innovative solution methodology based on relaxation and Logic-Based Benders Decomposition (LBBD) enhanced with logic-based inequalities and a customized neighbourhood search through combinatorial Benders cuts. The proposed models and solution methodology provide a comprehensive decision-support framework for optimizing evacuation logistics, with a particular emphasis on addressing the diverse medical needs of the affected population. To demonstrate the benefits of the proposed approach, we compare its performance to that of *three* alternative heuristic strategies for locating shelters, determining fleet size, and making routing decisions. The first strategy, called Closest-First (CF), evacuates the closest assembly areas to the vehicles first, transferring evacuees to the closest suitable medical centre. The second strategy, referred to as Shortest-Time-Window-First (STWF), gives priority to assembly areas with shorter time windows. The third strategy, called Most-Crowded-First (MCF), prioritizes the most crowded assembly areas.

To the best of our knowledge, the present paper is the first to offer a TSS optimization approach and a solution methodology for a setting that integrates a capacitated facility location problem and a stochastic multi-vehicle routing problems with hard time windows in the context of wildfire evacuation. This integration represents a significant advancement in evacuation modelling, being designed maximize evacuation within limited time windows in a cost-efficient manner while accounting for a range of critical factors and sources of uncertainty. These include road availability that may hinder evacuation routes, unknown count and distribution of affected population classified according to varying injury severity levels and medical needs, capacity constraints, and multiple types of vehicles each with its own technical specifications and operational constraints.

The remainder of this paper is organized in six sections. Section 2 reviews the relevant literature. Section 3 introduces the particular problem to be modelled and solved, including our main assumptions. Section 4 formulates the TSS optimization approach. Section 5 described the proposed solution methodology. Section 6 reports on extensive numerical experiments and provides descriptions of the three benchmark policies. Finally, Section 7 presents some concluding remarks and outlines potential directions for future research.

## 2. Literature Review

In general, the actions taken against the threat of natural disasters such as wildfires are categorized into three broad groups: pre-disaster preparation; disaster response; and post-disaster recovery (Zhang et al. 2021). In the context of wildfires, the literature on pre-disaster preparation is mostly centred around damage reduction, preventive actions and risk mitigation topics such as fuel break planning (Weise et al. 2023, Clark et al. 2023), prediction of fire spread and early fire detection using advance technologies (Carta et al. 2023, Zaidi 2023, Nur et al. 2023), and simulation of possible behavioural trends and assessment of evacuation protocols (Zhao et al. 2022).

On what regards disaster response, a substantial body of literature concerns fire control and the development of firefighting strategies. Readers are referred to the comprehensive literature reviews by Dhall et al. (2020) and Kamyabniya et al. (2024) for additional details. However, evacuation planning, especially for wildland-urban interface (WUI) communities, has always been one of the most important problems in this area. Various models have been developed to simulate wildfires and the evacuation process. For example, WUI-NITY (National Fire Protection Association 2024) and WUIVAC (Dennison et al. 2007).

Important decisions in the context of wildfire evacuation are whether or not and when to evacuate. Cova et al. (2011) developed an optimization-based decision-making tool to determine when and under what condition evacuation should be recommended to residents of a fire-affected zone as opposed to other options such as shelter-in-place. Extensive research has been done on stay-or-go decisions as well. Readers are referred to the comprehensive review by McLennan et al. (2019) for further information and references.

To support evacuation-related decisions, simulation and optimization models, individually or combined, have been used in the literature for purposes such as resource allocation, route selection, facility location and minimizing potential damage. Zhou & Erdogan (2019), for example, proposed an optimization approach that aims to minimize the number of individuals at risk in areas at high-risk of wildfire and thus reduce the overall costs associated with fire containment and property damage. The model optimizes emergency resource allocation decisions under capacity constraints while considering potential wildfire spread scenarios influenced by varying wind conditions. In another study, Gan et al. (2016) proposed a method for evacuation planning that integrates optimization with traffic simulation to determine the optimal evacuation timing and routes for evacuees. Through case studies based on data from wildfire-prone areas in Victoria, Australia, the authors demonstrated the benefits of their approach with respect to existing ad-hoc evacuation strategies.

Goerigk et al. (2014) developed an integrated model to simultaneously determining the location of shelters, bus routes for public transport and individual traffic routes while optimizing evacuation time. To efficiently solve the model within reasonable computational times, the authors developed a genetic algorithm, offering a heuristic approach to manage the complexity of the problem. Zhao et al. (2015) describe a multi-objective model designed to assign residents to earthquake shelters using the Beijing's Chaoyang district, China, as a case study. The model seeks to minimize the total weighted evacuation time from residential areas to designated shelters, while including constraints regarding shelter capacity and service radius, for three different scenarios concerning the number of evacuees. To solve the model, the authors adapted a particle swarm optimization algorithm, introducing a von Neumann structure in early iterations and a global structure in later stages. Although their work is not related to wildfire evacuation, and it does not account for time windows, it emphasizes the importance of minimizing evacuation time as a core objective.

A common aspect of the above-mentioned studies is that all assume residents are capable of self-evacuation and thus limited attention is given to supported-evacuation, which concerns vulnerable people who cannot independently evacuate the affected area (Flores et al. 2020, 2023). Flores et al. (2023) studied the challenges faced by individuals unable to evacuate independently due to physical, emotional, or logistical barriers. The authors modelled the evacuation process as a dynamic network of safe areas (shelters and hospitals) and unsafe areas (pick-up points), incorporating factors like health priority, dynamic arrivals, and heterogeneous vehicles. Using a mixed-integer goal programming model, their approach sought to maximize the number of evacuees, particularly high-priority ones, while minimizing evacuation time, costs and unmet supply needs. The model was validated through a case study based on the Saddleridge fire in California that took place in 2019.

Although the model in Flores et al. (2023) is comprehensive and includes both evacuation planning and logistics decisions for the required supplies, it lacks incorporating uncertainty concerning the number of evacuees, travel times, availability of routes and service capacity. Furthermore, the temporary shelters used for accommodating low-priority evacuees are considered already built, thus excluding facility location decisions from the problem.

Closely related to our study are the works by Shahparvari et al. Shahparvari et al. (2017), Shahparvari & Abbasi (2017) and Kamyabniya Kamyabniya (2022). Shahparvari et al. developed capacitated vehicle routing models to support last-minute wildfire evacuations within time windows, using heuristic approaches and comparing them to genetic algorithms across experiments and a real-world case—the 2009 Black Saturday wildfire. Their model assumes pre-established shelters and does not explicitly account for fire-impacted route availability, though reliability scores are assigned to routes. In a follow-up study, they incorporated route disruptions via a robust optimization framework. Kamyabniya, on the other hand, proposed a mixed-integer programming model for supported evacuation, jointly optimizing temporary shelter locations and vehicle routing. His model, formulated as a time-step-based two-stage stochastic problem, was ultimately solved using meta-heuristic techniques.

Our study distinguishes itself from these works through both its modeling framework and solution methodology. Unlike the robust approaches adopted by Shahparvari et al. and Kamyabniya, we employ a two-stage stochastic programming model that explicitly captures uncertainty through multiple scenarios.

We also propose an exact decomposition-based solution method, as opposed to the meta-heuristics used in the aforementioned studies. Additionally, our model incorporates a higher level of operational realism by including features such as split pickups and deliveries, continuous-time vehicle operations without return-to-depot constraints, vehicle-specific loading/unloading times, and first-stage fleet sizing decisions—factors that are either absent or handled differently in previous studies.

To better understand the proposed solution methodology, it is necessary to mention the foundational papers on Benders Decomposition (BD). BD was initially introduced by Benders (1962) and later extended by Geoffrion (1972) to cover a broader range of problems. It has since become one of the most widely adopted algorithms for solving large-scale optimization problems including complicating integer and binary variables. When these variables are fixed, the resulting model becomes amenable to decomposition into smaller sub-problems that are significantly easier to solve. A substantial body of literature has been devoted to developing variants of this algorithm to either accelerate its convergence to optimality or extend its applicability. For a comprehensive review, refer to Rahmaniani et al. (2017).

One of the primary challenges with BD is when it is applied to mixed-integer programming (MIP) problems involving both continuous and discrete decision variables in the sub-problems (SPs). While it has been demonstrated that for specific problems valid optimality and feasibility cuts and cuts derived from the dual counterpart of the LP-relaxation can be effective, a generalizable framework for generating useful cuts for this class of problems has yet to be established. To address this challenge, Codato & Fischetti (2006) proposed a variant of BD known as Combinatorial BD (CBD), which generates cuts based on the structure of the Restricted Master Problem (RMP) rather than on duality theory. CBD was originally developed for MIP models with binary variables in the RMP and feasibility SPs with no objective function. A couple of years earlier, Hooker & Ottosson (2003) introduced Logic-Based BD (LBBD), which incorporates problem-specific knowledge into cut generation to iteratively reduce the solution space of the RMP.

In this paper, given the structure of the problem, we decompose each proposed optimization model into a binary RMP and multiple MIP SPs. We adopt ideas from both CBD and LBBD, introducing a customized neighbourhood search that employs combinatorial cuts. Additionally, to further narrow the solution space and expedite the convergence of the proposed algorithm, we incorporate logic-based inequalities and a warm start.

## 3. Problem Definition

To plan for evacuation, we assume that there are $|I|$ known locations referred to as assembly areas inside the fire zone. The number and the priority of patients waiting in each assembly area for rescue are unknown. Outside the fire zone, there are $|J|$ medical facilities. Medical facilities and assembly areas, collectively called nodes, are connected to each other by routes. Among the medical facilities, $|J_1|$ are already established hospitals that can serve patients of all priority levels ($h \in \{1\text{: high}, 2\text{: low}\}$) and $|J_2|$ are potential locations considered for shelter setup, such as schools, malls, and stadiums. Shelters are suitable for serving low-priority evacuees only. The capacity of medical facilities are limited for each priority level and can vary according to the scenario. Although the potential locations for shelters are known, the optimization model finds the optimal number and selected locations of shelters. We consider multiple types of ground and aerial vehicles. $V_1$ is the set of ambulances, $V_2$ buses, and $V_3$ helicopters. Each vehicle type has a different speed, capacity, and activation cost. Helicopters can directly travel between locations and only ambulances are eligible to transport high-priority patients as they typically need medical attention and/or vital devices on their way to medical facilities. We assume known upper bounds on the number of vehicles of each type. Each vehicle is allowed to have a maximum number of stops $Max_T$ and vehicle-related variables in each stop are indexed by $t \in \{1, 2, ..., Max_T\}$. Evacuation vehicles are allowed to pickup patients from multiple assembly areas (split pickup) and deliver patients to multiple medical facilities (split delivery).

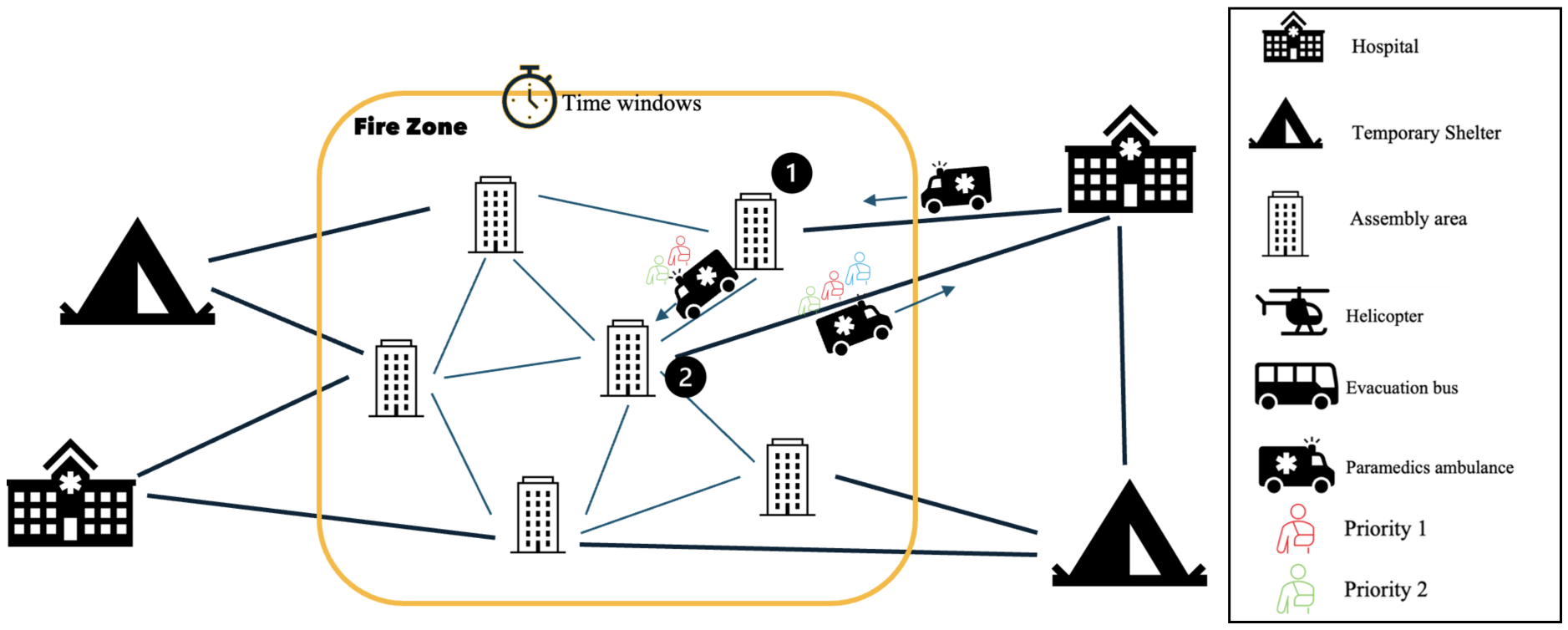

Figure 1: A sample wildfire evacuation network

The problem setting can be represented as a network where assembly areas, hospitals and shelters are the nodes, and the roads connecting them to each other are the links. Roads have unlimited capacity but their availability can be affected by the fire. To incorporate the different sources of uncertainty into our models, we define a finite set of $|S|$ scenarios. All stochastic parameters and decisions variables related to these scenarios are indexed by $s \in S$. The first model to be introduced later assumes fixed travel times of $T^{vs}_{ll'}$ for vehicle $v$ between location $l$ and $l'$ but it may vary depending on scenario $s$. The second model to be introduced considers non-stationary travel times reflecting the reality that roads might become partially disrupted or congested during the evacuation leading to longer travel times. Depending on the origin of the fire, wind direction, and vegetation density, each assembly area $i$ is subject to a time window of $TW^s_i$ after which it cannot be visited. Figure 1 demonstrates a sample network diagram depicting the aforementioned components as well as an example of split pickup.

Under scenario $s$, an uncertain number of $Pa^s_{ih}$ patients with priority $h$ are waiting in assembly area $i$. Each evacuation vehicle $v$ is initially stationed at a specific medical facility $j$ under scenario $s$ (determined by a binary parameter $IL^s_{vj}$). To enhance the realism of the models, in addition to travel time, for each vehicle $v$, extra time of $LV_{vi}$ and $UV_{vj}$ may be required for loading patients at assembly area $i$ or unloading them at medical facility $j$, respectively. For ground vehicles, travelling between two locations is only possible if there is a direct connection between them (determined by a binary parameter $APA_{ll'}$). The model also accounts for various costs and penalties: these include the cost of setting up a shelter at a given location ($OC_j$), the cost of activating and using a vehicle to transport patients ($OV_v$), penalties for failing to evacuate patients based on their priority level ($PN_h$), travelling unit cost ($\psi_v$), and penalties for delaying the evacuation beyond the acceptable time limits ($\xi$). Main decision variables used in the model are grouped into first-stage and second-stage variables. $X_v$ and $Y_j$ determine the activation of a vehicle and a shelter, respectively. Routing of vehicles across the network is determined by $VT^{vt}_{ll't}$. The number of pickups and drop-offs are determined by $GI^{vt}_{his}$ and $GO^{vt}_{hjs}$, respectively, and the arrival times are captured by $AT^{ls}_{vt}$. Detailed description of the notation used in the proposed models in this paper are set out in Table 1.

Table 1: Mathematical Notation

| Set | Description | Index | Description | Subset | Description |
|---|---|---|---|---|---|
| $I$ | Assembly areas | $i$ | An assembly area, $i \in I$ | $J_1$ | Hospitals, $J_1 \subseteq J$ |
| $J$ | Medical facilities | $j$ | A hospital or shelter, $j \in J$ | $J_2$ | Shelters, $J_2 \subseteq J$ |
| $L$ | Locations, $L = I \cup J$ | $l$ | A location, $l \in L$ | $V_G$ | Ground vehicles, $V_G = V_1 \cup V_2$ |
| $H$ | Priority levels = $\{1, 2\}$ | $h$ | A priority level, $h \in H$ | $V_k$ | Vehicles of type $k$ |
| $K$ | Vehicle types = $\{1, 2, 3\}$ | $k$ | A vehicle type, $k \in K$ | $T_v$ | Stops of vehicle $v$ |
| $V$ | Vehicles | $v$ | A vehicle, $v \in V$ | | |
| $T$ | Vehicle stops | $t$ | A vehicle stop, $t \in T$ | | |
| $S$ | Scenarios | $s$ | A scenario, $s \in S$ | | |

| Parameter | Description |
|---|---|
| $\pi^s$ | Probability of occurrence of scenario $s$. |
| $TW^s_i$ | Time window for the evacuation at assembly area $i$. |
| $\alpha^s$ | The percentage by which travel time on a road that is partially disrupted or congested increases under scenario $s$. |
| $Pa^s_{ih}$ | Number of patients at assembly area $i$ from level $h$ under scenario $s$. |
| $T^{vs}_{ll'}$ | Travel time between $l$ and $l'$ under scenario $s$ for vehicle $v$. |
| $DT^s_{ll'}$ | Direct travelling time between $l$ and $l'$ for helicopter under scenario $s$. |
| $M_1$ | Maximum number of shelters. |
| $M_2$ | Maximum number of vehicles of type $V_1$, $V_2$, Helicopter that we have. |
| $OC_j$ | Operational cost for setting up a shelter in location $j \in J_2$. |
| $OV_v$ | Operational cost of activating vehicle $v$. |
| $CV_v$ | Capacity of vehicle $v$. |
| $C^s_{jh}$ | Capacity of hospital or shelter $j$ for accepting patients from severity level $h$. |
| $PN_h$ | Penalty of leaving a patient with severity $h$ behind. |
| $APA_{ll'}$ | 1 if location $l$ and $l'$ are immediately connected to each other. |
| $DT_{ll'}$ | Direct travelling time between $l$ and $l'$ for helicopter. |
| $IL^s_{vj}$ | 1 if the initial location of vehicle $v$ is $j \in J$ under scenario $s$. |
| $LV_{vi}$ | Loading time of vehicle $v$ in location $i$. |
| $UV_{vj}$ | Unloading time of vehicle $v$ in location $j$. |
| $RDT^s_{ll'}$ | Time at which the route connecting locations $l$ and $l'$ becomes partially disrupted under scenario $s$, increasing its travel time. |
| $\xi$ | A penalty term to minimize the arrival times. |
| $\psi_v$ | Travelling unit cost for vehicle type $v$. |
| $M$ | A very large value. |
| $Max_T$ | Maximum number of stops. |

| Variable | Description |
|---|---|
| $Y_j$ | (1st-stage) 1 if a shelter is set up in location $j$, 0 otherwise. |
| $X_v$ | (1st-stage) 1 if vehicle $v$ participates in the evacuation, 0 otherwise. |
| $VT^{vt}_{ll's}$ | (2nd-stage) 1 if vehicle $v$ travels between locations $l$ and $l'$ after its $t$-th stop under scenario $s$, 0 otherwise. |
| $GI^{vt}_{his}$ | (2nd-stage) No. of priority-$h$ patients picked up in assembly area $i$ at the $t$-th stop of vehicle $v$ under scenario $s$ ($GI^{vt}_{his} \in Z^+_0$). |
| $GO^{vt}_{hjs}$ | (2nd-stage) No. of priority-$h$ patients dropped off at facility $j$ at the $t$-th stop of vehicle $v$ under scenario $s$ ($GO^{vt}_{hjs} \in Z^+_0$). |
| $AT^{ls}_{vt}$ | (2nd-stage) Time of vehicle $v$'s $t$-th stop at location $l$ under scenario $s$ ($AT^{ls}_{vt} \geq 0$). |
| $RS^{vt}_{ll's}$ | (2nd-stage) 1 if vehicle $v$ takes partially disrupted route from location $l$ to $l'$ after its $t$-th stop under scenario $s$, 0 otherwise. |
| $PV^{vt}_{hs}$ | (2nd-stage) No. of priority-$h$ patients on vehicle $v$ after its $t$-th stop under the scenario $s$ ($PV^{vt}_{hs} \in Z^+_0$). |
| $PU^{is}_{vt}$ | (2nd-stage) 1 if vehicle $v$ picks up patients from assembly area $i$ during its $t$-th stop under scenario $s$, 0 otherwise. |
| $DO^{js}_{vt}$ | (2nd-stage) 1 if vehicle $v$ drops off patients at medical facility $j$ during its $t$-th stop under scenario $s$, 0 otherwise. |
| $Q^{vt}_{ls}$ | (2nd-stage) 1 if stop $t$ at location $l$ is the final stop of vehicle $v$ under scenario $s$, 0 otherwise. |

## 4. Mathematical Formulation

To formulate this problem, we propose a two-stage stochastic (TSS) optimization approach. The first stage consists of decisions regarding the activation of vehicles and location of shelters. In the second stage, the model optimizes evacuation routing decisions depending on the realization of the stochastic parameters. The first model seeks to minimize objective function (1), which includes the costs associated with vehicle activation and shelter setup as well as the penalties associated with unevacuated individuals, travelling costs and unnecessary delays.

$$\textbf{Model 1:}\ \text{Min } Z = \sum_{j\in J_2} OC_j Y_j + \sum_v OV_v X_v + \sum_s \pi^s \left[ \sum_i \sum_h PN_h \left( Pa^s_{ih} - \sum_v \sum_t GI^{vt}_{his} \right) + \sum_l \sum_{l'} \sum_v \sum_t \psi_v T^{vs}_{ll'} VT^{vt}_{ll's} + \xi \sum_v \sum_t \sum_l AT^{ls}_{vt} \right] \tag{1}$$

subject to constraints (2) to (32) below.

$$\sum_v \sum_t GO^{vt}_{his} \leq 0 \qquad \forall h, i, s \tag{2}$$

$$\sum_v \sum_t GI^{vt}_{hjs} \leq 0 \qquad \forall h, j, s \tag{3}$$

$$\sum_v \sum_t GO^{vt}_{h=1,js} \leq 0 \qquad \forall j \in J_2, s \tag{4}$$

Constraints (2) and (3) ensure that pickups and drop-offs are done at the right locations. Constraint (4) ensures that priority level 1 patients are not delivered to shelters.

$$\sum_v \sum_t GO^{vt}_{hjs} \leq C^s_{jh} \qquad \forall h, j \in J_1, s \tag{5}$$

$$\sum_v \sum_t GO^{vt}_{hjs} \leq C^s_{jh} \times Y_j \qquad \forall h, j \in J_2, s \tag{6}$$

Constraints (5) and (6) ensure that the capacities of the medical facilities are respected.

$$\sum_{l'} VT^{v,t=0}_{ll's} \leq IL^s_{vl} \times X_v \qquad \forall v, l, s \tag{7}$$

$$NV^{v,t=0}_{ls} = IL^s_{vl} \qquad \forall v, l, s \tag{8}$$

$$VT^{vt}_{ll's} \leq APA_{ll'} \qquad \forall v \in V_G, l, l', t, s \tag{9}$$

$$NV^{vt}_{ls} = \sum_{l'} VT^{v,t-1}_{l'ls} \qquad \forall v, l, t \in T_v, s \tag{10}$$

$$\sum_l NV^{vt}_{ls} \leq 1 \qquad \forall v, t, s \tag{11}$$

Constraints (7) and (8) set the initial locations of the evacuation vehicles. Constraints (9) to (11) ensure that vehicles only travel on existing links and visit only one node per stop.

$$GI_{his}^{v,t=0} = 0 \qquad \forall v, h, i, s \tag{12}$$

$$GO_{hjs}^{v,t=0} = 0 \qquad \forall v, h, j, s \tag{13}$$

$$GI_{hls}^{vt} + GO_{hls}^{vt} \leq CV_v \times NV_{ls}^{vt} \qquad \forall v, h, l, t \in T_v, s \tag{14}$$

Constraints (12) to (14) ensure that the evacuation vehicles start their trips empty and that the total number of patients on each vehicle must not exceed its capacity.

$$\sum_h GI_{his}^{vt} \leq CV_v \times PU_{vt}^{is} \qquad \forall v, i, t, s \tag{15}$$

$$PU_{vt}^{is} \leq \sum_h GI_{his}^{vt} \qquad \forall v, i, t, s \tag{16}$$

$$\sum_h GO_{hjs}^{vt} \leq CV_v \times DO_{vt}^{js} \qquad \forall v, j, t, s \tag{17}$$

$$DO_{vt}^{js} \leq \sum_h GO_{hjs}^{vt} \qquad \forall v, j, t, s \tag{18}$$

$$\sum_t \sum_i GI_{his}^{vt} = \sum_t \sum_j GO_{hjs}^{vt} \qquad \forall v, h, s \tag{19}$$

Constraints (15) to (18) set binary variables for indicating pickups and drop-offs, which will be useful later on in the calculation of arrival times. Constraint (18) states that, for each vehicle across all its stops, the total number of pickups and drop-offs must be equal.

$$PV_{hs}^{v,t=0} = 0 \qquad \forall v, h, s \tag{20}$$

$$\sum_h PV_{vs}^{ht} \leq CV_v \qquad \forall v, t, s \tag{21}$$

$$PV_{hs}^{vt} = PV_{hs}^{v,t-1} + \sum_l \left(GI_{hls}^{vt} - GO_{hls}^{vt}\right) \qquad \forall h, v, t \in T_v, s \tag{22}$$

Constraints (20) to (22) determine how the number of patients on each vehicle is updated after each stop.

$$\sum_v \sum_t GI_{his}^{vt} \leq Pa_{ih}^{s} \qquad \forall i, h, s \tag{23}$$

Constraint (23) limits the total number of pickups to the number of patients waiting in each assembly area.

$$\sum_l VT_{ll's}^{vt} = \sum_{l''} VT_{l'l''s}^{v,t+1} + Q_{l's}^{v,t+1} \qquad \forall v, l', t \in T_v \setminus \{Max_T\}, s \tag{24}$$

$$Q_{is}^{vt} = 0 \qquad \forall v, i, t, s \tag{25}$$

Constraints (24) and (25) determine that if a vehicle is activated, it should end its trip eventually.

$$AT^{ls}_{v,t=0} = 0 \qquad \forall v, l, s \tag{26}$$

$$AT^{ls}_{vt} \leq M \times NV^{vt}_{ls} \qquad \forall v, l, t, s \tag{27}$$

Constraints (26) and (27) impose that arrival times start at 0 and are set to zero for unvisited nodes.

$$AT^{is}_{vt} \leq \left(TW^{s}_{i} - LV_{iv}\right) \times NV^{vt}_{ls} \qquad \forall v, i, t, s \tag{28}$$

Constraint (28) forbids travelling to nodes whose time windows are closed.

$$\begin{aligned} AT^{ls}_{vt} \geq & \sum_{l'} AT^{l's}_{v,t-1} + \sum_{l'} LV_{l'v} PU^{l's}_{v,t-1} + \sum_{l'} UV_{l'v} DO^{l's}_{v,t-1} + \sum_{l'} T^{vs}_{l'l} VT^{v,t-1}_{l'ls} \\ & - M\left(1 - \sum_{l'} VT^{v,t-1}_{l'ls}\right) \quad \forall v \in V_G, l, t \in T_v, s \end{aligned} \tag{29}$$

$$\begin{aligned} AT^{ls}_{vt} \geq & \sum_{l'} AT^{l's}_{v,t-1} + \sum_{l'} LV_{l'v} PU^{l's}_{v,t-1} + \sum_{l'} UV_{l'v} DO^{l's}_{v,t-1} + \sum_{l'} DT^{s}_{l'l} VT^{v,t-1}_{l'ls} \\ & - M\left(1 - \sum_{l'} VT^{v,t-1}_{l'ls}\right) \quad \forall v \in V_3, l, t \in T_v, s \end{aligned} \tag{30}$$

Constraints (29) to (30) determine the arrival time of vehicles at different locations.

$$GI^{vt}_{h=1,ls} = 0 \qquad \forall v \notin V_1, l, t \in T_v, s \tag{31}$$

$$GO^{vt}_{h=1,ls} = 0 \qquad \forall v \notin V_1, l, t \in T_v, s \tag{32}$$

Constraints (31) and (32) make sure that priority-1 patients are only transported using ambulances.

Since the activation and routing of vehicles as well as the capacitated facility location of shelters are both strongly NP-hard problems, using off-the-shelf software to solve Model 1 at once can be very time-consuming for reasonably-sized problems. To address this issue, we propose the hybrid solution methodology based on Benders Decomposition described next.

## 5. A Hybrid Benders Decomposition Solution Approach

We propose a hybrid solution methodology combining ideas from relaxation, Benders Decomposition (BD) and neighbourhood search as illustrated in Figure 2. The main issue with adopting the classic form of BD for our models is that the second-stage variables are not all continuous–they consist of a mix of continuous, binary and integer variables–making the sub-problems MIPs as well. Consequently, deriving feasibility and optimality cuts is not straightforward, as their derivation typically relies on dual formulations, which are not directly available for mixed-integer programming models. The integer variant of BD offered by Laporte & Louveaux (1993) is not applicable to all MIP models as it was originally designed to address discrete sub-problems with a particular structure, such as containing only binary variables. To address this challenge, we propose a hybrid solution methodology built on ideas from Combinatorial BD (CBD) and Logic-Based BD (LBBD) further enhanced by the use of logic-based inequalities and a customized neighbourhood search.

### 5.1. logic-based Inequalities

To limit the solution space of the first-stage variables and expedite convergence to the optimal solution, we introduce several logic-based inequalities that remove some of the sub-optimal solutions ensuring that the optimal solution is not cut out.

- *Symmetry Removal.* Since same-type vehicles are considered identical, Constraint (33) orders the activation of vehicles of the same type according to their index.

$$X_{v+1} \leq X_v \qquad \forall v \in \{1, \ldots, |V_k| - 1\}, k \in \{1, 2, 3\} \tag{33}$$

- *Lower Bounds on the Number of Vehicles Required.* Let $MP_h$ be the maximum total number of priority $h$ patients waiting to be evacuated across all scenarios and $H_k$ be the set of priority levels that each vehicle of type $k \in \{1\text{: ambulance, } 2\text{: bus, } 3\text{: helicopter}\}$ can serve. Assuming that $CV_k$ is the capacity of a vehicle of type $k$ and that an upper bound $AV_k$ on the number of vehicles of type $k$ is available, then we can impose a lower bound on the number of vehicles of type $k$ as follows:

$$NV_k = \begin{cases} \min\left\{AV_k, \left\lceil \frac{\sum_{h \in H_k} MP_h}{CV_k \times \frac{Max_T}{2}} \right\rceil \right\} & \text{if } k \in \{1, 2\} \\ \min\left\{AV_k, \left\lceil \frac{RP}{CV_k \times \frac{Max_T}{2}} \right\rceil \right\} & \text{if } k = \{3\} \text{ and } \sum_{k' \in K \setminus \{3\}} NV_{k'} > \sum_{k \in K \setminus \{3\}} AV_{k'} \\ 0 & \text{if } k = \{3\} \text{ and } \sum_{k' \in K \setminus \{3\}} NV_{k'} \leq \sum_{k \in K \setminus \{3\}} AV_{k'} \end{cases} \tag{34}$$

where $Max_T$ is the maximum number of stops a ground vehicle can make, $\lceil \cdot \rceil$ denotes the ceiling function, which rounds a number up to the nearest integer, and $RP$ is the maximum total number of patients, across all scenarios, that ground vehicles are not expected to evacuate.

$$RP = \max\left\{0, \sum_{k \in K \setminus \{3\}} \left( \sum_{h \in H_k} MP_h - NV_k \times CV_k \times \frac{Max_T}{2} \right) \right\}.$$

**Theorem 1.** *$NV_k$ is a valid lower bound on the number of vehicles of type $k \in \{1, 2\}$ for full evacuation.*

The proof of Theorem 1 is provided in the appendix.

Given that the activation cost for helicopters is much higher than that for ground vehicles, a reasonable lower bound on the number of helicopters can be determined when there are insufficient ground vehicles to evacuate everyone.

Thus, Constraint (35) makes sure that solutions with fewer vehicle activations than the estimated minium are not explored.

$$NV_k \leq \sum_{v \in V_k} X_v \leq AV_k \qquad \forall k \tag{35}$$

- *Lower Bound on the Number of Shelters Required.* Since shelters can only receive low-priority evacuees, it is reasonable to assume that some shelters will be set up if the capacity of hospitals is insufficient to serve all low-priority patients. Consequently, Constraint (36) guarantees the activation of a minimum number of shelters to provide sufficient added capacity to serve low-priority patients. Additionally, Constraint (37) limits the number of shelters to be less that a pre-defined maximum upper bound ($USH$).

$$\sum_{j \in J_2} C^s_{j,h=2} Y_j \geq MP_{h=2} - \sum_{j \in J_1} C^s_{j,h=2} \quad \forall s \in S \tag{36}$$

$$\sum_{j \in J_2} Y_j \leq USH \tag{37}$$

- *Shelter Location Priority.* If additional capacity to serve low-priority patients is needed, or if having more shelters is convenient to reduce travel times, priority can be given to setting up shelters based on a weighting scheme. A weight $WY_j$ can be assigned to shelter $j$ according to its travelling distance (based off of evacuation *buses*) from assembly areas depending on the number of evacuees there. As a result, shelters that are closer to assembly areas hosting a larger number of low-priority evacuees are more attractive choices. This prioritization assumes that shelters are similar in terms of activation cost.

$$WY_j = \sum_{i,s} Pa^s_{i,h=2} \times \frac{1}{T^{vs}_{ij}} \qquad \forall j \in J_2, v \in V_2 \tag{38}$$

The potential shelter locations can then be sorted in ascending order based on the values of $WY_j$. If we assume, without loss of generality, that the elements of set $J_2$ are already sorted, then Constraint (39) is imposed.

$$Y_j \leq Y_{j-1} \qquad \forall j \in \{2, \ldots, |J_2|\} \tag{39}$$

The above logic-based inequalities are added to the RMP to avoid unnecessary exploration.

### 5.2. Relaxation-Based Sequential Decomposition

Under certain conditions, vehicles of the same type can be considered identical and independent. Specifically, if all vehicles of a given type share the same capacity, speed, activation cost, traveling cost, and initial location, treating them as identical seems reasonable. Additionally, if the $\xi$ penalty in the objective function is sufficiently small—serving only to ensure that the model minimizes the arrival times in Constraints (29) and (30)—then vehicle trips remain independent and can be interchangeable. Although all vehicles operate within the same network, picking up evacuees from shared assembly areas and competing for capacities from the same hospitals and shelters, optimizing each vehicle's itinerary individually leads to the optimization of the entire fleet. As a result, vehicle "cooperation" (i.e., making the routing decisions of the entire fleet at the same time from an optimization point of view) does not provide additional benefits in this setting. We will further analyze the strength of this assumption and identify the conditions under which it may not hold. Below, we explain how the decomposition of Model 1 into restricted master problem and sub-problems is done.

**Restricted Master Problem (RMP):**

$$\min Z_{\text{RMP}} = \sum_{j \in J_2} OC_j Y_j + \sum_{v} OV_v X_v \tag{40}$$

subject to Constraints (33), (35), (36), (37), and (39).

This brings us to a feasible solution for the first-stage variables that can be iteratively refined by searching its neighbourhood with a radius of $\mathcal{R}$. This neighbourhood is called "Trust Region".

**Sub-Problems (SPs):**

SPs determine the optimal routing for each vehicle and are solved through a relaxation-based, sequential decomposition approach. In Model 1, most of the constraints are specified to each vehicle except for constraints (5), (6), and (23). These three constraints are linking constraints that, if relaxed, make the problem easily decomposable by vehicle. However, relaxing these constraints leads to an infeasible solution with respect to the original model. To handle this issue, the following constraints are individually imposed for each sub-problem:

$$\sum_{t} GO_{hjs}^{vt} \leq C_{jh}^{s} \qquad \forall h, v, s, j \in J_1 \tag{41}$$

$$\sum_{t} GO_{h=2,js}^{vt} \leq C_{jh}^{s} Y_j \qquad \forall v, s, j \in J_2 \tag{42}$$

$$\sum_{t} GI_{his}^{vt} \leq Pa_{ih}^{s} \qquad \forall h, v, s, i \tag{43}$$

Therefore, the SP for vehicle $v$ is formulated as follows:

$$\begin{aligned} \min \theta_v^s = \pi^s \Bigg[ \sum_{i} \sum_{h} PN_h \Bigg( \min \Big\{ Pa_{ih}^s, \Big\lceil CV_k \times \frac{Max_T}{2} \Big\rceil \Big\} - \sum_{t} GI_{his}^{vt} \Bigg) \\ + \psi_v \sum_{l} \sum_{l'} \sum_{t} T_{vs}^{ll'} VT_{ll's}^{vt} + \xi \sum_{t} \sum_{l} AT_{vt}^{ls} \Bigg] \end{aligned} \tag{44}$$

subject to Constraints (2) to (32) for vehicle $v$, excluding Constraints (5), (6), and (23), plus Constraints (41), (42), and (43). In the objective function (44), $\lceil CV_k \times \frac{Max_T}{2} \rceil$ is the theoretical maximum number of patients vehicle $v$ can pickup. The reason for including the minimum between this term and $Pa_{ih}^s$ in the objective function is to avoid penalizing for individuals who were not supposed to evacuate by the current vehicle. For example, when there are 18 priority-1 patients in the fire zone, and the first activated ambulance can pickup at most 10 patients in total, leaving the remaining for the second ambulance, the sub-problem for ambulance 1 would show a huge unreasonable penalty for not rescuing that 8 patients who are supposed to be picked up by ambulance 2. Since the penalty term $PN_h$ is very large, inclusion of only $Pa^s_{ih}$ in the objective function would divert the focus of the model away from other terms in the objective function related to reducing travel times and unnecessary delays.

### *5.3. Neighbourhood Search through Combinatorial Cuts*

In this sub-section, an algorithm is described that conducts an $\mathcal{R}$-radius neighbourhood search of the most recent solution of the first-stage variables. In each iteration $r$ of the algorithm, according to a perturbation rule, we alter the value of each of first-stage variables by a maximum of $s$ units using the combinatorial cuts below and make the current solution temporarily infeasible to force the model to explore:

$$\textbf{Vehicle Perturbation:} \quad \sum_{v \in S^r} (1 - X^r_v) + \sum_{v \notin S^r} X^r_v \leq \mathcal{R} \tag{45}$$

$$\textbf{Shelter Perturbation:} \quad \sum_{j \in S^r} \left(1 - Y^r_j\right) + \sum_{j \notin S^r} Y^r_j \leq \mathcal{R} \tag{46}$$

$$\textbf{Temporary Cut (TC):} \quad \sum_{v \in S^r} (1 - X^r_v) + \sum_{v \notin S^r} X^r_v + \sum_{j \in S^r} \left(1 - Y^r_j\right) + \sum_{j \notin S^r} Y^r_j \geq 1 \tag{47}$$

In this expression, the set $S^r = (\bar{X}^r, \bar{Y}^r)$ is the set of first-stage variables having non-zero values (i.e., activated vehicles and shelters). TC ensures the current solution is temporarily infeasible, forcing the model to explore new solutions.

**Perturbation Rule through logic-based Benders:** In each iteration, there are four possible actions: adding a shelter, removing a shelter, activating a vehicle, and deactivating a vehicle. If in the current solution there are unevacuated individuals, then depending on their priority level, the lower bound for the corresponding vehicle type defined in expression (34) is increased in the RMP to make sure everyone is evacuated. Constraints (45) to (47) iteratively assist the model to generate new solutions. The model explores $NS^+$ neighbours, which involves changing the value of inactive variables from zero to one, thereby activating additional shelters or vehicles. Constraints (33) and (39) define the activation order for vehicles and shelters, guaranteeing unique solutions. During the search, some solutions are dominated by a better solution. To avoid re-generating these sub-optimal solutions, they are stored in a set $BS^r$ for iteration $r$ of the search algorithm. Constraint (48) is imposed on the RMP for iterations $r > 1$.

$$\sum_{v \in BS^r} (1 - X^r_v) + \sum_{v \notin BS^r} X^r_v + \sum_{j \in BS^r} \left(1 - Y^r_j\right) + \sum_{j \notin BS^r} Y^r_j \geq 1 \tag{48}$$

If activation/deactivation of more vehicle/shelter leads to a worse objective function, any activation/deactivation in that specific direction are considered sub-optimal and cut out from the solution space. In other words, upper bounds defined in Constraints (34) and (37) are updated after each bad move (e.g., $Z_r > Z_{r-1}$ where $Z_r = Z_{RMP} + \sum_{s,v} \theta^s_v$ is the overall objective of the problem in iteration $r$). Similarly, after each good move (e.g., $Z_r < Z_{r-1}$) the corresponding lower bounds on the first-stage variables are updated. For example, when adding a new shelter leads to a significant decrease in the travel times and generates a better overall objective function, the lower bound on the number of shelters must be updated to avoid unnecessary exploration. The search continues until all the possible solutions within the updated lower/upper bounds are visited. In other words, until the RMP gets infeasible due to Constraint (47), which means the best solution so far is the global optimal for the fist-stage variables. Thus, the magnitude of search radius $\mathcal{R}$ does not affect the execution time or the number of iterations as long as $s$ is large enough. The entire solution approach explained above is summarized in Figure 2. An algorithmic description of the solution approach is also presented in the appendix.

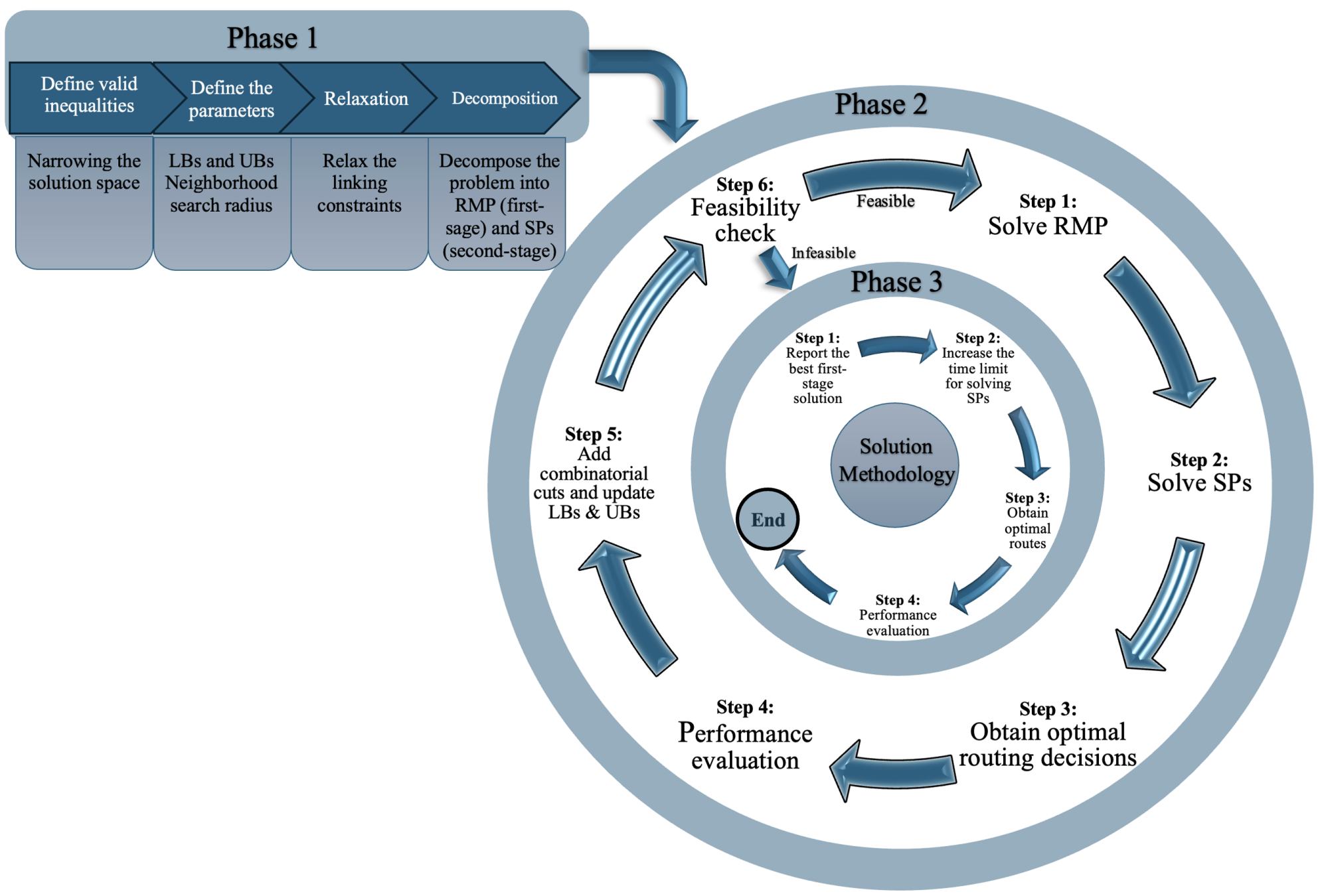

Figure 2: Hybrid Solution Methodology

An advantage of using the proposed solution approach is that, in the intermediate iterations, the SPs do not need to be solved to optimality, allowing us to set short time limits or optimality gap thresholds to expedite the solution progress. Once the best solution for the $X$ and $Y$ variables are found, SPs can be re-solved, allowing more time to each SP to further refine the final solution.

We alternatively explored the use of Lagrangian relaxation to tackle our optimization model. However, the results were not desirable. As mentioned earlier, Constraints (5), (6), and (23) are linking constraints that after being relaxed make the problem decomposable by vehicle. According to Lagrangian relaxation, we can penalize the violation of the aforementioned constraints in the objective function. However, this form of relaxation is not useful because Constraint (23) is the same as the third term in the objective function (1), where the number of unevacuated individuals is penalized. Therefore, assuming that $\lambda_{ih}^{s}$ is the Lagrangian multiplier corresponding to this constraint, the following term will appear in the Lagrangian objective function:

$$\sum_{i}\sum_{h}(PN_h - \lambda_{ih}^{s})\left(Pa_{ih}^{s} - \sum_{v}\sum_{t} GI_{his}^{vt}\right) \tag{49}$$

Thus, no matter what kind of sub-gradient algorithm is used to obtain the optimal multipliers, $(PN_h - \lambda_{ih}^{s})$ will take the solution to two extremes. If the penalty for leaving patients behind is less than that of violating the constraint, then the model will not allow vehicles to pickup patients from there. On the other hand, if the penalty for violating the constraint is smaller, then more patients than actually exist will be picked up from that location, leading to an infeasible solution.

### *5.4. A key assumption and model extension*

The solution methodology proposed in this section is based on two layers of decomposition (e.g., by scenarios and vehicles). The decomposition by vehicles has been proposed based on the assumption that vehicles of the same type are identical, and no benefit is obtain from their cooperation, and thus they can be treated independently. In this subsection, we propose a version of the model where this assumption no longer holds. We introduce uncertainty associated with the availability of routes connecting different nodes of the evacuation network, and the travel times on the affected routes. In a real world setting, it is reasonable to assume that during an evacuation, travel times are affected by traffic congestion and fire spread. In the latter case, routes might become partially or completely disrupted, leading us to consider time windows for them similar to assembly areas.

To incorporate these changes, we included secondary travel times as well as route time windows into the model. We assume each route has two time windows. The first is a soft time window, $RDT^{s}_{ll'}$, after which the route between $l$ to $l'$ becomes partially disrupted under scenario $s$, increasing the route travel time by $\alpha$%. The second is a hard time window, $TWR^{s}_{ll'}$, after which the route connecting $l$ to $l'$ is reached by the fire and becomes completely unavailable. To model these, we introduce a new binary variable, $RS^{vt}_{ll's}$, which is equal to 1 if vehicle $v$ takes the route from location $l$ to $l'$ after its $t^{th}$ stop under partial disruption in scenario $s$. We also add a new dimension to the arrival time variable, $AT^{ll's}_{vt}$, which is the origin of the trip to have control over the routes vehicles take to reach their destinations in each trip and watch for partial and/or complete route disruptions.

With the above-mentioned changes, constraints (26) to (30) are updated to constraints (50) to (55):

$$AT^{ll's}_{(v,t=0)} = 0 \qquad \forall s, v, l, l' \tag{50}$$

$$\sum_{l'} AT^{l'ls}_{vt} \leq M \times NV^{vt}_{ls} \qquad \forall s, v, t, l \tag{51}$$

Constraints (50) and (51) state that arrival times start at 0 and are set to 0 for unvisited nodes.

$$\sum_{l} AT^{lis}_{vt} \leq (TW^{s}_{i} - LV_{iv})NV^{vt}_{ls} \qquad \forall s, v, t, i \tag{52}$$

$$\sum_{l''} AT^{l''ls}_{vt} \leq TWR^{s}_{ll'} + M(1 - VT^{vt}_{ll's}) \qquad \forall s, v, t, l, l' \tag{53}$$

Constraints (52) and (53), respectively, forbid traveling to nodes and using routes whose time windows are closed.

$$AT_{vt}^{ll's} \geq \sum_{l''} AT_{v,t-1}^{l''ls} + PU_{v,t-1}^{ls} LV_{lv} + DO_{v,t-1}^{ls} UV_{lv} + T_{ll'}^{s} \left(1 + \alpha RS_{ll's}^{v,t-1}\right)$$

$$- M \times \left(1 - VT_{ll's}^{v,t-1}\right) \quad \forall s, l, l', v \in V_G, t \in T_v \tag{54}$$

$$AT_{vt}^{ll's} \geq \sum_{l''} AT_{v,t-1}^{l''ls} + PU_{v,t-1}^{ls} LV_{lv} + DO_{v,t-1}^{ls} UV_{lv} + DT_{ll'}^{s}$$

$$- M \times \left(1 - VT_{ll's}^{v,t-1}\right) \quad \forall s, l, l', v \in V_3, t \in T_v \tag{55}$$

Constraints (54) to (55) determine the arrival time of vehicles at different locations.

To incorporate the surge in travel times due to partial route disruptions, we impose the following logical condition as linear constraints: *"If the arrival time of vehicle $v$ at location $l$ in its $t^{th}$ stop is greater than $RDT_{ll'}^{s}$, then a binary variable $RS_{vt}^{ll's}$ in addition to the main movement variable, $VT_{ll's}^{vt}$, must be equal to 1"*. Constraints (56) and (57) state that:

$$\sum_{l''} AT_{vt}^{l''ls} - RDT_{ll'}^{s} \leq M \times RS_{vt}^{ll's} + M \times (1 - VT_{ll's}^{vt}) \qquad \forall s, v, t, l, l' \tag{56}$$

$$VT_{ll's}^{vt} \geq RS_{vt}^{ll's} \qquad \forall s, v, t, l, l' \tag{57}$$

The binary variable $RS_{vt}^{ll's}$ will help us to account for the extra time vehicles spend due to partial road disruptions or congestion. Thus, the updated objective function is given by (58).

$$\textbf{Model 2:} \ \min Z = \sum_{j \in J_2} OC_j Y_j + \sum_{v} OV_v X_v + \sum_{s} \pi^s \left[ \sum_{i} \sum_{h} PN_h \left( Pa_{ih}^{s} - \sum_{v} \sum_{t} GI_{his}^{vt} \right) \right.$$

$$\left. + \sum_{l} \sum_{l'} \sum_{v} \sum_{t} \psi_v T_{ll'}^{vs} \left(VT_{ll's}^{vt} + \alpha RS_{ll's}^{vt}\right) + \xi \sum_{v} \sum_{t} \sum_{l} \sum_{l'} AT_{vt}^{ll's} \right] \tag{58}$$

One secondary travel times are incorporated into the formulation, the assumption of vehicle independence no longer holds, as vehicles of the same type can coordinate to balance workload and minimize travel times on partially disrupted routes. For example, a solution where vehicle 1 completes 20 trips and vehicle 2 completes 8, which may appear optimal under sequential decomposition, could be suboptimal. A better solution might involve vehicle 1 stopping after 14 trips to reduce travel on congested routes, allowing vehicle 2 to take 6 additional trips before soft time windows close. Despite of this, scenario decomposition remains useful through combinatorial Benders decomposition, with neighborhood search iteratively refining first-stage decisions. This adjustment increases sub-problem (SP) complexity, as each SP now involves multiple vehicles. Nonetheless, applying phases 2 and 3 of our solution approach illustrated in Figure 2 reduces the complexity of the full formulation, especially for instances involving a large number of scenarios. However, scaling to larger networks with a large number of vehicles remains a significant challenge.

## 6. Numerical Experiments

In this section, we provide details of numerical experiments conducted to discuss the practical benefits and showcase the applicability of the proposed models and solution methodology. All the numerical experiments were conducted on a 2022 MacBook Pro with an M1 Pro chip and 16 GB of memory.

### *6.1. Benchmark Policies*

To evaluate the performance of the proposed approach, we identified three alternative evacuation policies based on common sense and current practices. The three alternative policies explained here share the same strategy for the first-stage decisions, however, routing strategies are different. Particularely, reaching the assembly areas are based on the following logic:

1. ***C**losest assembly area **F**irst (CF):* Shortens travel times by choosing the closest assembly area for each trip with respect to the current location of the vehicle.

2. *Assembly area with the **S**hortest **T**ime-**W**indow **F**irst (STWF):* Minimizes the risk of exceeding time-windows by choosing assembly areas with the shortest available time.

3. ***M**ost **C**rowded assembly area **F**irst (MCF):* Balances the evacuation by choosing assembly areas with the larger remaining population of evacuees at any point of time.

Once the vehicles are at assembly areas, under all three benchmark heuristics, they are instructed to transfer the patients they picked up to the closest medical facilities provided that it has sufficient capacity for the specific priority level of the patients on the vehicle.

Vehicles are activated based on the number of trips they can make and their used capacity in each trip. Earlier, this was determined using parameter $Max_T$ that was set as an input parameter in the optimization model due to the requirements of the formulation. If $Max_T$ is known, the upper bound on the number of trips for each vehicle is $\lceil \frac{Max_T}{2} \rceil$. Otherwise, a reasonable estimate of the number of trips, as a starting point, can be obtained based on the time windows and the average travel times between assembly areas and medical facilities:

$$\text{Estimated number of trips for vehicle type } k = \left\lceil \frac{1}{2} \times \frac{\frac{1}{|S||I|} \sum_{s\in S} \sum_{i\in I} TW_i^s}{\frac{1}{|S||I||J|} \sum_{s\in S} \sum_{i\in I} \sum_{j\in J} T_{i,j}^{ks}} \right\rceil \tag{59}$$

Based on our analysis, working with average travel time provides reasonable results, while incorporating the worst case (i.e., longest travel times and/or shortest time windows) would only cause waste of resources (i.e., unnecessary activation of vehicles). Thus, the proposed vehicle activation rule considers the population of the priority level(s) that this type of vehicle can evacuate and the estimated number of trips.

$$\text{Estimated number of active vehicles of type } k = \left\lceil \frac{\sum_{h\in H_k} \sum_{i\in I} \sum_{s\in S} Pa_{i,h}^s}{|S| \cdot CV_k \cdot \text{Number of trips}_k} \right\rceil \tag{60}$$

where $H_k$ is the set of priority level(s) a vehicle of type $k$ can evacuate. The above equation is based on the capacity of the vehicle, therefore provides a lower bound on the number of vehicles required. It assumes that vehicles are loaded to their capacity at each assembly area, which is not always the case as patients are typically more scattered in the final moments of the evacuation. To deal with this issue, we allow the number of vehicles to be increased iteratively if the number determined by the Equation (60) is not sufficient to evacuate everyone under the evacuation policy in use.

Locating shelters depends on two important factors: whether a shelter must be open due to the lack of service capacity for low-priority evacuees; and whether opening a shelter leads to a significant decrease in the travel times between assembly areas and medical facilities. If the first scenario applies, the minimum number of shelters necessary to satisfy Constraint (36) is activated. Otherwise, the set of locations that contributes the most to shortening the average travel times between assembly areas and medical facilities, while considering the activation cost, is selected. Assuming that there are $|J_2|$ potential locations for shelters, let $Q$ be the set of all $2^{|J_2|}$ possible location combinations with $Q_i$ be the $i^{th}$ combination. Let $\mathcal{T}$ be the corresponding average travel times between shelters and active medical facilities and $\mathcal{T}_i$ be for the $i^{th}$ combination. Finally, let $Q_0$ and $\mathcal{T}_0$ represent the scenario in which no shelter is activated. Equation (61) determines the best combination of shelters to open with respect to the above-mentioned criteria.

$$\arg\min_{i \in Q} \left( \sum_j OC_j |Q_i| + \psi_{v \in V_2} (\mathcal{T}_i - \mathcal{T}_0) \right) \tag{61}$$

In the following subsections, we describe multiple numerical experiments to demonstrate the applicability of the proposed models and solution methodology.

### 6.2. *A Small Instance*

We designed a small instance of the problem consisting of three assembly areas (A1, A2, A3), two hospitals outside the fire zone (H4, H5), and two potential shelter locations within a closer distance to the assembly areas compared to the hospitals (S6, S7). Table 2 provides information on the number and capacity of vehicles of each type as well as the count and distribution of the patients with two priority levels under three distinct scenarios. Vehicles are allowed to make a maximum of 20 stops. As discussed earlier, the maximum number of stops allowed, $Max_T$, is required in the formulation as an input. However, provided its value is set large enough, it does not impact the optimal solution. Table 2 also provides the time windows associated with each assembly area.

### 6.3. *Results*

The results for this instance suggest opening one shelter in location S6 due to both the lack of low-priority medical capacity and to expedite the evacuation by reducing the travel times. This solution is consistent with the shelter activation rule of the benchmark policies. Figure 3a illustrates the trade-off between shelter setup costs and the potential reduction in average travel times from assembly areas to all the medical facilities. As can be seen, although building two shelters at S6 and S7 reduces the average travel times, considering the shelter activation cost, the most cost-efficient option is to set up one shelter at S6.

Table 2: Characterization of a small instance of the problem

| **Vehicle Type** | **Count** | **Capacity** | | |
|---|---|---|---|---|
| Ambulance | 4 | 1 | | |
| Bus | 4 | 6 | | |
| Helicopter | 2 | 2 | | |
| **Affected population** | | **Assembly area** | | |
| **Location** | **A1** | **A2** | **A3** | |
| High-priority | {2, 5, 4}* | {5, 4, 4} | {3, 5, 2} | |
| Low-priority | {15, 15, 17} | {15, 19, 22} | {18, 24, 23} | |
| Time-window [min] | {838, 846, 745} | {678, 600, 798} | {741, 827, 711} | |
| **Capacity** | | **Medical facility** | | |
| **Location** | **H4** | **H5** | **S6** | **S7** |
| High-priority | {19, 17, 19} | {15, 13, 13} | — | — |
| Low-priority | {9, 12, 15} | {10, 14, 13} | {60, 63, 62} | {59, 59, 63} |

*: {Scenarios 1, Scenario 2, Scenario 3}

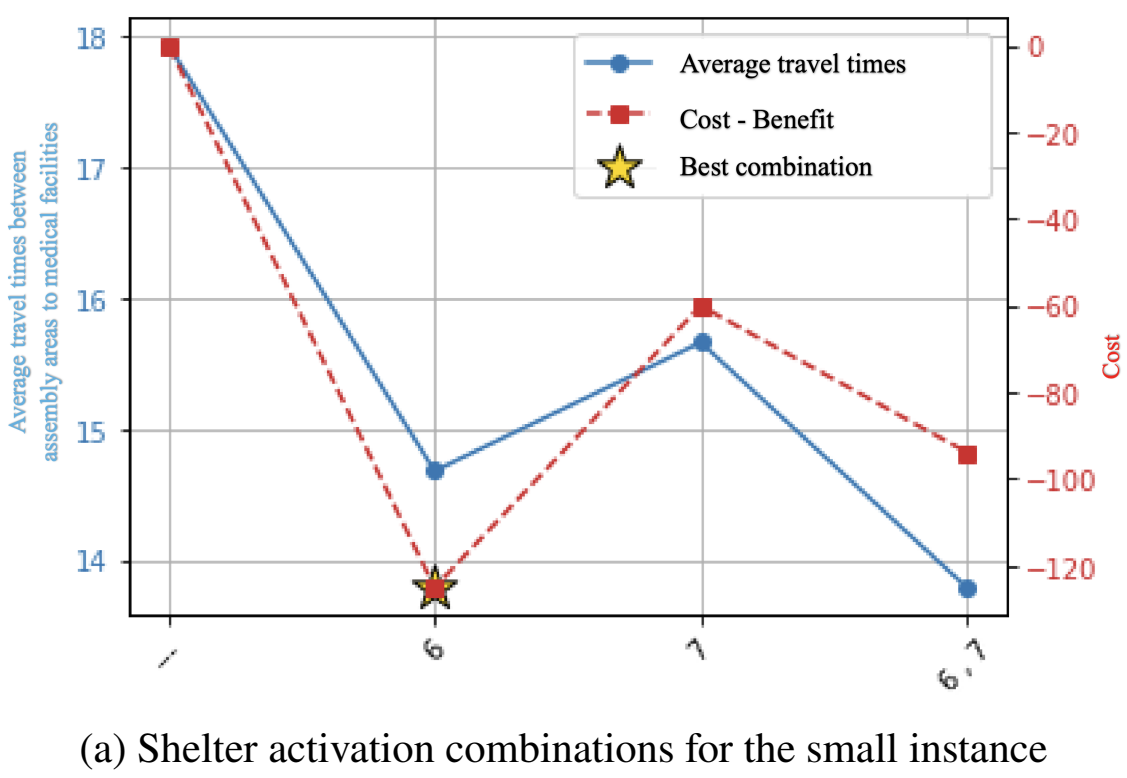

(a) Shelter activation combinations for the small instance

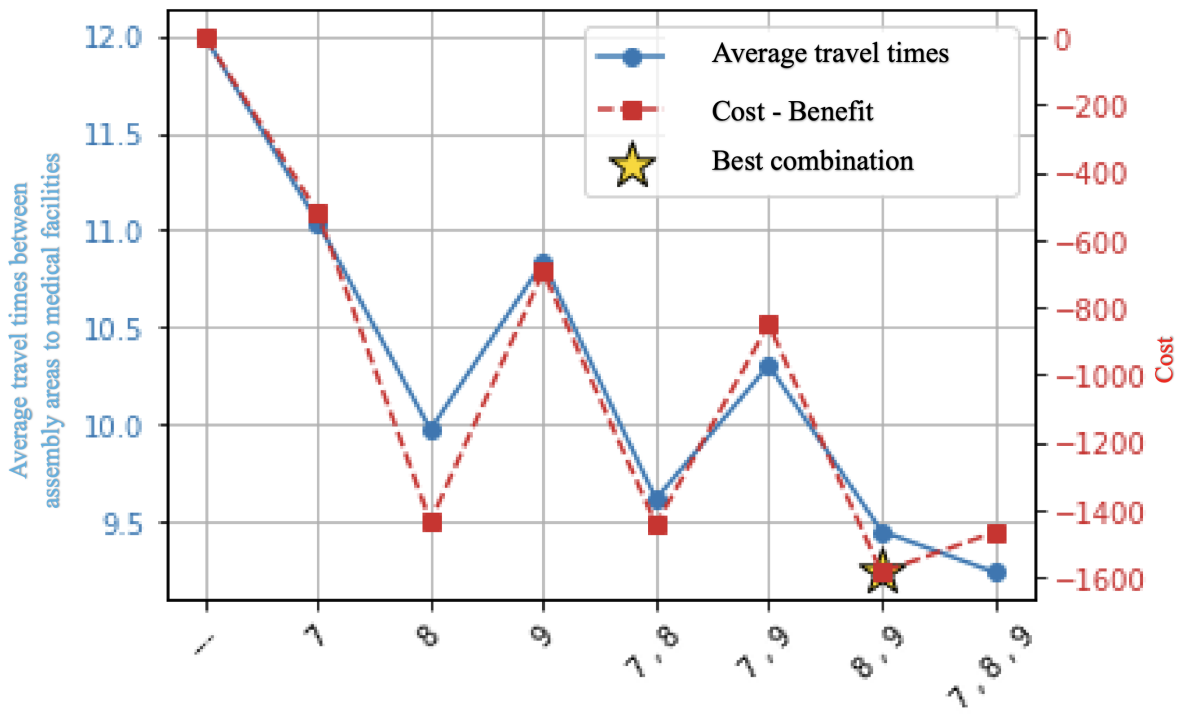

(b) Shelter activation combinations for the practical instance

Figure 3: Combined figure showing two different shelter activation analyses.

In terms of vehicle activation, the model activates two ambulances, two buses, and zero helicopters. Table 3 shows the itinerary of one of the buses under Scenario 1. The bus starts its journey from location H4, which is a hospital, and throughout its trips, it transfers patients from assembly areas A1, A2, and A3 to shelter S6, and ends its trip there. It is worth noting that the solution demonstrated in Table 3 is obtained directly using the commercial solver (i.e., no decomposition) with a time limit of 5 hours. After 5 hours, the aforementioned solution had a 23% optimality gap. To find a better solution faster, we implemented a warm start. Specifically, we initialized the solver with the best solution obtained using the CF heuristic, which performed better than the other alternatives for this instance. With a warm start, the optimality gap was reduced to 11.2% in 5 hours. However, the optimality gap is mostly related to a poor lower bound rather than to the poor quality of the solution. This issue will be further examined later when applying the proposed solution methodology to more complex instances. Additionally, we will present further results demonstrating the benefits of using a warm start for both the default solver and the proposed solution approach.

Figure 4 illustrates how Model 1 optimizes the use of resources and completes the evacuation efficiently. As can be seen, for this instance, the solution obtained from solving Model 1 outperforms the three heuristics with respect to the evacuation time (bottom right subgraph) as well as the average travelling distance taken by vehicles (left subgraph). The figure also shows how the four activated vehicles are used across different scenario (top right subgraph).

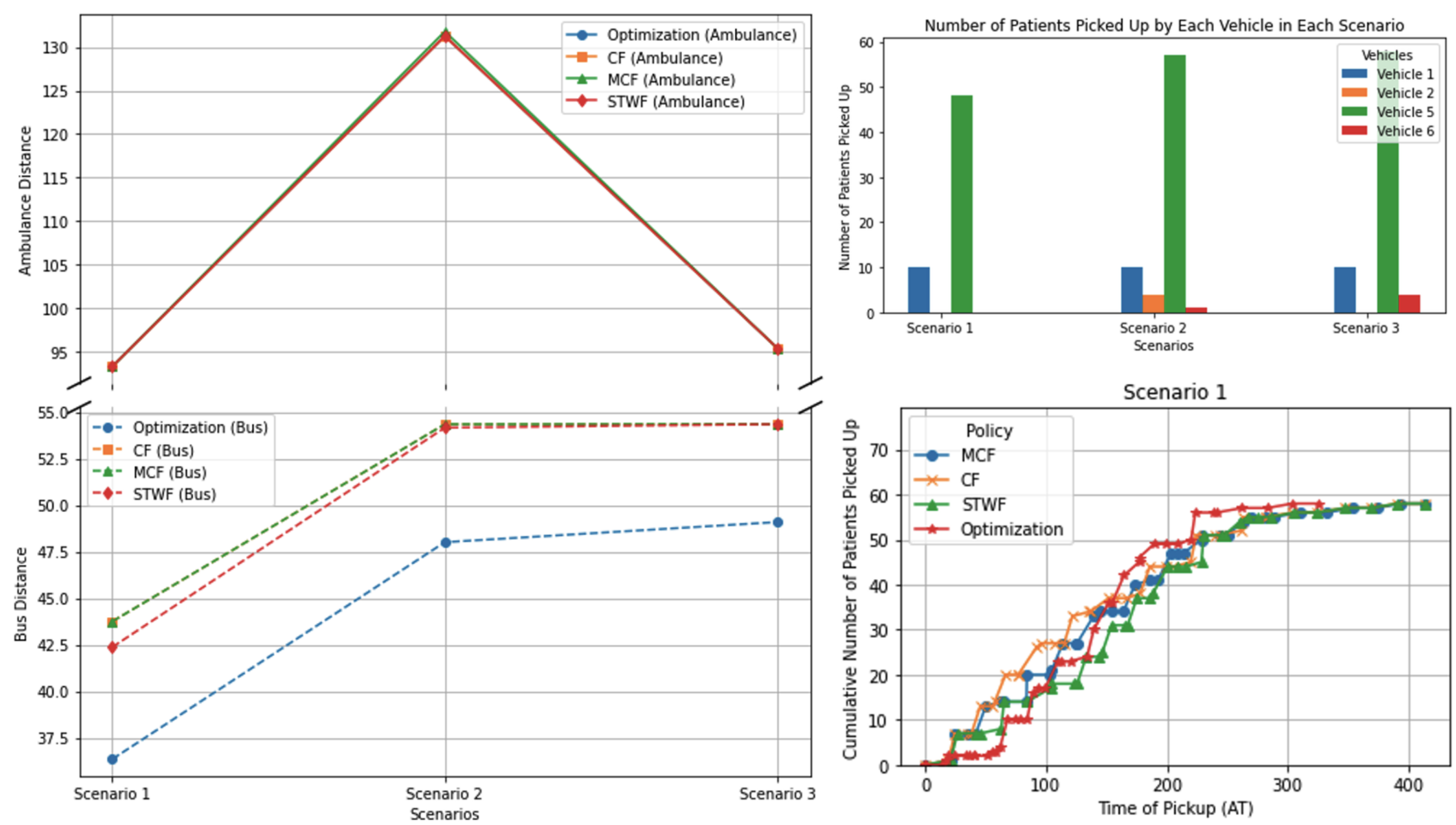

Figure 4: Results obtained for the small instance

Table 3: Itinerary of vehicle 5 (a bus) under Scenario 1 for the small instance

| Stop | 1 | 2 | 3 | 4 | 5 | 6 | 7 | 8 | 9 | 10 | 11 | 12 | 13 | 14 | 15 | 16 | 17 | 18 |
|---|---|---|---|---|---|---|---|---|---|---|---|---|---|---|---|---|---|---|
| Location | H5 | A2 | A1 | S6 | A3 | S6 | A3 | S6 | A3 | S6 | A1 | S6 | A1 | S6 | A2 | S6 | A2 | S6 |
| Pickups | 0 | 3 | 3 | 0 | 6 | 0 | 6 | 0 | 6 | 0 | 6 | 0 | 6 | 0 | 6 | 0 | 6 | 0 |
| Drop-offs | 0 | 0 | 0 | 6 | 0 | 6 | 0 | 6 | 0 | 6 | 0 | 6 | 0 | 6 | 0 | 6 | 0 | 6 |
| Arrival Time [min] | 0 | 23 | 35 | 50 | 61 | 71 | 82 | 92 | 102 | 113 | 128 | 143 | 158 | 173 | 192 | 211 | 230 | 250 |

### *6.4. Larger Instances*

The primary challenge with Model 1, as previously noted, was the long solving times and large optimality gaps. Even small instances could not be solved to optimality within a reasonable time, and as problem size increases, not only does optimality become unattainable, but solutions obtained after extensive computation still exhibit large optimality gaps.

To evaluate the effectiveness of the proposed solution approach, we test four randomly generated configurations of increasing size and compare the performance of the proposed methodology against those of the direct solution and alternative evacuation policies in terms of their objective value and solving time.

Table 4 highlights the importance of warm-starting the solver with the best heuristic solution. Without a warm start, the quality the direct solution remains worse than that of the heuristic solutions, which are obtained in seconds. However, a warm-start enables the solver to refine the solution. Our results show that when first-stage variables align between all solution approaches, the improvement is modest (e.g., Configurations 1 and 2, with gains of 3.2% and 5.9%, respectively). The improvement primarily stems from the allowance of split pickups in the optimization approach, which mostly benefits multi-passenger vehicles (e.g., buses) in the final trips of the last activated vehicle.

When the first-stage decisions obtained by the alternative heuristics are suboptimal, the gain of using the optimization approach is more considerable (e.g., Configurations 3 and 4). This indicates that the second-stage decisions obtained by the alternative heuristics are effective, and the advantage of the optimization approach mostly due to better first-stage decisions. This is however subject to whether we can assume vehicles of the same type operate independently. Additionally, our methodology achieves in minutes what for the default solver takes 24 hours to compute. Another key advantage of the proposed approach is that the subproblems (SPs) in our approach are significantly simpler than the original model, allowing most SPs to be solved to optimality. This suggests that the overall optimality gap is primarily due to weak lower bounds rather than poor solution quality.

Overall, it is worth noting that even for the proposed solution approach the computation time for large instances can still be high as in each iteration of the algorithm, there are $|V| \times |S|$ SPs each of which involves solving an NP-hard vehicle routing problem with time windows (VRPTW). Nonetheless, the proposed methodology provides significant computational advantages as summarized below:

- The logic-based inequalities reduce the solution space significantly and provide managerial insights regarding the minimum levels of resources required.
- Solving the entire problem directly is equivalent to solving a multi-VRPTW indexed by $|S|$ scenarios, which is much harder compared to multiple single VRPTWs that can be solved in parallel for $|S|$ scenarios and sequentially for $|V|$ vehicles.
- Given the lower complexity of the SPs compared to the entire formulation, most of them are solved either to optimality or with a minimal optimality gap, providing explanation for the large optimality gap of the default solver (direct solution).
- Given that the entire solution process is broken into multiple steps, the proposed methodology can allocate more time to solve SPs and thus to refine and improve their solutions once the value of the first-stage variables are known.

When both soft and hard time windows are imposed on routes, vehicles of same type can no longer be treated independently. Optimizing the routing of vehicles considering their simultaneous operation (vehicle cooperation) can cause benefits such as balancing workload and minimizing commutes on routes that are partially affected by the fire. This collaboration also reduces the overall evacuation time by ensuring a more balanced distribution of trips, preventing scenarios where some vehicles make significantly more trips than others. Methodologically, we adapt the proposed solution approach by eliminating the sequential vehicle decomposition stage, decomposing the original problem into just $|S|$ SPs instead, and refining first-stage decisions through combinatorial cuts in a neighborhood search framework.

Table 4 demonstrates that in this case, the optimization approach outperforms the benchmark heuristics to a greater extent. Given the added complexity of Model 2 compared to Model 1, solving it directly using the default solver is significantly harder. Without a warm start, the solver fails to provide useful solutions, and even with one, refinement is slow and solving times remain high. While the added complexity also affects the SPs in our solution approach, they remain easier to solve, allowing us to obtain high-quality solutions in much less time.

Table 4: Comparison of the performance of Model 1 and Model 2 under four different network configurations.

| Config. | $\lvert I\rvert, \lvert J_1\rvert, \lvert J_2\rvert, \lvert V_1\rvert, \lvert V_2\rvert, \lvert V_3\rvert$ | Model 1 | | | | | Model 2 | | | | |
|---|---|---|---|---|---|---|---|---|---|---|---|
| | | Default Solver (24 hours) | | Benders Decomposition | | Best Heuristic | Default Solver (24 hours) | | Benders Decomposition | | Best Heuristic |
| | | With Warm Start | Without Warm Start | **Best Obj.** | **Time (s)** | **(Obj.)** | With Warm Start | Without Warm Start | **Best Obj.** | **Time (s)** | **(Obj.)** |
| 1 | 3, 2, 2, 4, 3, 1 | 2525 (11%) | 2889 (23%) | 2525 | 121 | STWF (2609) | 2638 (23%) | 3238 (38%) | 2516 | 393 | STWF (2774) |
| 2 | 4, 3, 2, 5, 4, 1 | 2836 (29%) | 3357 (39.4%) | 2830 | 168 | CF (3017) | 2987 (36%) | 3469 (44%) | 2870 | 1438 | CF (3276) |
| 3 | 5, 3, 3, 5, 4, 2 | 3357 (28%) | 4387 (45%) | 3368 | 523 | CF (3952) | 3628 (40%) | 4147 (47%) | 3429 | 2890 | CF (4317) |
| 4 | 6, 4, 4, 6, 5, 3 | 3469 (13%) | 4154 (28%) | 3501 | 882 | CF (3798) | 3721 (27%) | 3822 (29%) | 3683 | 3776 | CF (3983) |

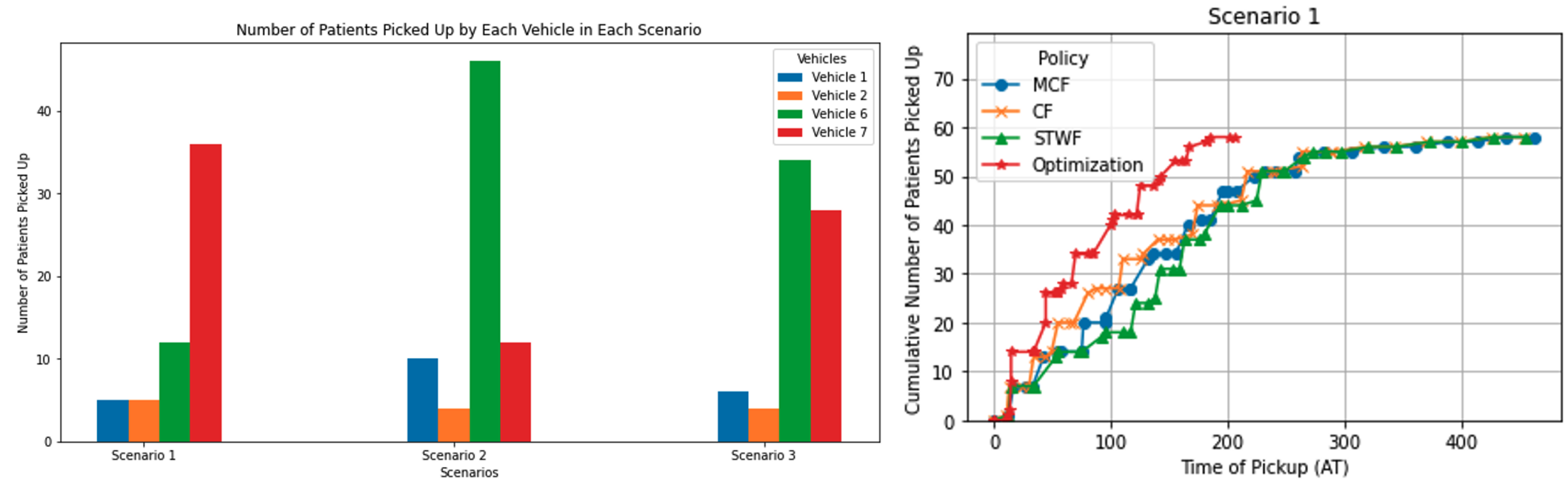

Figure 5: Performance of Model 2 in reducing evacuation time and balancing the workload

Figure 5 further illustrates how Model 2 handles the small instance introduced in Subsection 6.2. Comparing with Figure 4, it is evident that the vehicles of the same type (e.g., vehicles 1 and 2, and vehicles 6 and 7) share their workload more evenly. Moreover, at any given moment, the optimization approach evacuated more people than the heuristics. This highlights Model 2's ability to balance vehicle workload effectively, leading to significantly shorter evacuation times.

### *6.5. Case Study - Roxborough Park*

In this section, we consider a case study based on the data associated with a large-scale fire evacuation drill presented by Gwynne et al. (2023) at Roxborough Park, Colorado, US a WUI community of almost 900 homes at potential risk of wildfires. In 2019, the authors collected data on evacuation times, route usage, and arrival times at a designated assembly areas. This data served as a critical resource for testing two evacuation models: WUI-NITY and the Evacuation Management System (EMS). The study underscored the significant impact of pre-evacuation delays and route selection on evacuation effectiveness, revealing that even with prior knowledge of the drill, some residents took over an hour to begin evacuating. A total of 143 households participated in the evacuation drill. The evacuation relied on the residents use of their own vehicles and adherence to evacuation protocols. In this section, we create a realistic supported evacuation problem based on this case study to evaluate the performance of our models.

### 6.5.1. Instance Creation

According to the 2022 US census, the population of Roxborough park is around 9,057 people. But there is no information available on the actual number of people who need evacuation support in case of wildfires. A quick search of the US census database suggests that close to 12.9% of the population are over the age of 65, some of whom may not be able to drive or use public transportation independently, and 5.8% of the population have a disability. Thus, it seemed reasonable to assume that approximately 10% of the population (nearly 900 people) would be in need of supported-evacuation. This vulnerable population may include individuals with disability, patients or injured people who need to be taken to a hospital and need medical attention on the way, elderlies, residents with no car, and tourists in the area.

We chose three locations, 1) a fire station in the east, 2) an intermediate school, and 3) a primary school, as potential assembly areas where residents can wait to be transferred to shelters or hospitals. Outside the fire zone, three locations (mostly schools) were chosen as potential places for temporary shelters. Also, three hospitals were chosen with a similar capacity of 80 beds to serve patients in need of medical attention (UCHealth Highlands Ranch Hospital, Sky Ridge Medical Center, and Advent Health Castle Rock). It is assumed within the hospitals or close to them there are places that can be used as already established shelters to accommodate low-priority evacuees. Figure 6 shows on the map of Roxborough park the aforementioned locations, and Table 5 shows the shortest distance between any two given locations according to Google Maps estimates from January 2025.

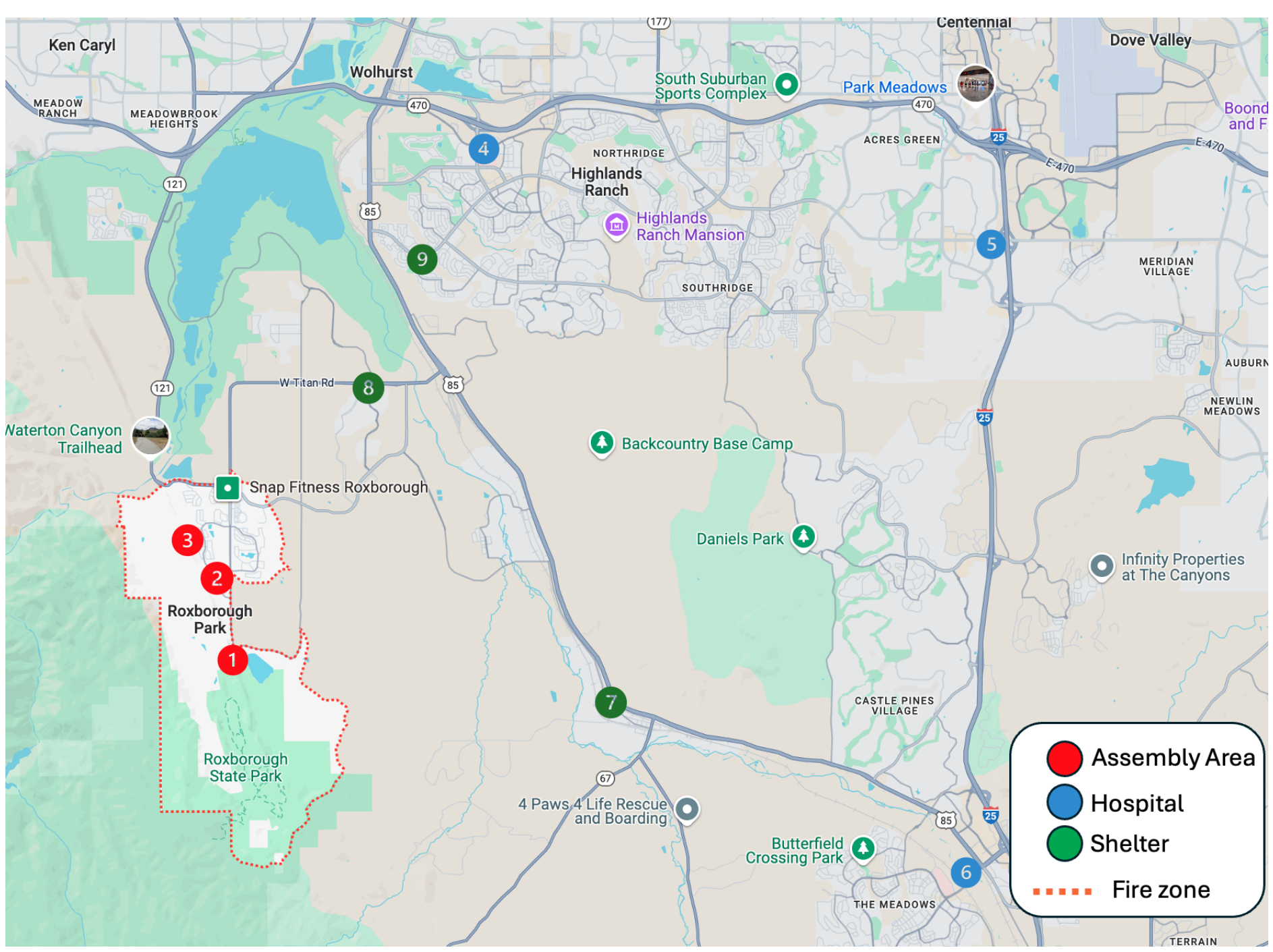

Figure 6: Map of Roxborough Park including the location of assembly areas, hospitals and shelters (Source: Google Maps)

We assumed the availability of 10 ambulances with capacity of one patient and average speed of 60 km/h, 20 buses with capacity of ten patients and average speed of 50 km/h, and 2 helicopters with capacity of two patients and average speed of 200 km/h.

Table 5: Distance between locations in kilometers. Letter A denotes an assembly area, letter H a hospital, and letter S a shelter.

| Origin | Destination | | | | | | | | |
|---|---|---|---|---|---|---|---|---|---|
| | A1 | A2 | A3 | H4 | H5 | H6 | S7 | S8 | S9 |
| A1: Fire Station | – | 2.4 | 3.9 | 18.5 | 35.5 | 28.4 | 19.4 | 9.3 | 19.1 |
| A2: Intermediate School | 2.5 | – | 1.5 | 17.1 | 34.7 | 26.7 | 18.0 | 9.7 | 15.4 |
| A3: Primary School | 3.0 | 1.8 | – | 19.1 | 33.8 | 26.1 | 17.1 | 7.0 | 14.5 |
| H4: UCHealth Highlands Ranch | 18.7 | 17.1 | 15.9 | – | 16.3 | 27.1 | 18.1 | 9.7 | 4.0 |
| H5: Sky Ridge Medical Center | 38.6 | 35.6 | 34.4 | 13.6 | – | 15.8 | 21.0 | 22.0 | 16.0 |
| H6: AdventHealth Castle Rock | 28.2 | 27.0 | 25.5 | 27.0 | 17.6 | – | 10.6 | 21.3 | 25.2 |
| S7: Sedalia Elementary School | 18.7 | 17.1 | 15.9 | 18.4 | 22.9 | 10.8 | – | 12.4 | 16.2 |
| S8: Primrose School | 9.6 | 8.1 | 6.9 | 9.5 | 22.0 | 21.3 | 12.6 | – | 7.8 |
| S9: Coyote Creek School | 19.3 | 15.3 | 14.1 | 3.8 | 16.0 | 24.8 | 16.2 | 7.8 | – |

Table 6: Details of the stochastic parameters used in instance creation.

| Parameter | Description | Stochastic Value |
|---|---|---|
| $TW_i^s$ | Time windows for assembly areas [min] | Uniform(400, 850) |
| $RDT_{ll'}^s$ | Time windows for routes [min] | Uniform(400, 850) |
| – | Fire origin | Random(south,west, north) |
| – | Proportion of high-priority patients | $Bi(P_1 = 0.05)$ |
| – | Population | $900 \times [1 + \text{Uniform}(-20\%, +20\%)]$ |
| $C_{jh}^s, j \in J_2$ | Shelter capacities | Uniform(400, 600) |
| $C_{j,h=1}^s, j \in J_1$ | Hospital capacities for high-priority evacuees | $80 \times \text{Uniform}(10\%, \%20)$ |
| $C_{j,h=2}^s, j \in J_1$ | Hospital capacities for low-priority evacuees | Uniform(80, 100) |

We randomly generate 10,000 scenarios corresponding to different realizations of the stochastic parameters according to the probability distributions in Table 6. In order to obtain a manageable number of more distinct scenarios, we classified these 10,000 scenarios according to four criteria: 1) fire origin, 2) average length of time-windows, 3) population level, and 4) proportion of high-priority patients. Fire origin has 3 possibilities (south, north, or west). Average length of the time windows is classified into 2 groups (long or short). Scenarios with a time window shorter than the mean of the corresponding probability distribution were labelled as ”short” and others as ”long”. For the two population-related criteria, we labelled scenarios with values within one standard deviation of the mean of the distribution as ”moderate” and anything below and above that as ”low” and ”high”, respectively. Thus these two criteria are categorized into 3 levels. In total, we have $2 \times 3 \times 3 \times 3 = 54$ different combinations defining our final set of aggregated scenarios (hereafter referred to as ”category”). For each category, we took the average value of the input parameter over all scenarios that fell into that particular category and round it to the nearest integer. The probability associated with each category was determined by dividing the total number of scenarios that fell into that category by 10,000.

*6.5.2. Results*

Solving Model 2, which is the more realistic formulation, provided a solution with 6 ambulances, 16 buses, and 0 helicopters. Among the three potential locations for opening shelters, two locations (S8 and S9), were chosen. This was also consistent with the shelter activation rule developed for the alternative benchmark heuristics. According to Figure 3b, opening three shelters indeed helps to reduce the average travel times between assembly areas and medical facilities. However, considering the total shelter activation cost together with the expected reduced average transportation cost that is provided by each alternative, opening two shelters at S8 and S9 is the most beneficial choice. The results shows that the full evacuation can be done within the time windows resulting in almost 10% improvement compared to benchmark approaches due to better routing decisions.

Graphical displays of the results related to this instance is presented in the appendix of this paper. According to the results, hospital H5 was less frequently used across categories as it is located at a farther distance. While hospital H4 was found more useful due to its closer distance to the fire zone. The two activated shelters in locations S8 and S9 were used to accept a large number of low-priority evacuees across all categories (more in categories with a larger affected population). As expected, the shelter located in location S8 was used more due to its closer distance to the fire zone. Figure 7 demonstrates that hospitals were mostly used for high-priority evacuees who need medical attention, except for hospital H4 that also accepted low-priority evacuees as it was located close to the fire zone.

In Figure 7, the five selected categories are structured so that each consecutive pair differs based on a single criterion used to classify the original 10,000 scenarios into 54 categories. First, consider categories 1 and 2, where the key distinction is in the proportion of high-priority patients. Category 1 has a higher proportion, leading to a slight increase in evacuation time as shown in the cumulative pickup curve. This is due to the greater reliance on ambulances, which have lower transport capacity than buses. Next, the comparison between categories 2 and 5 highlights the effect of overall population size. Category 5 has the smallest population size across all priority levels, resulting in lower utilization of shelters and hospitals and a lower upper bound in the cumulative pickup curve. Between categories 5 and 14, the key difference lies in the average length of the time windows. Category 5 has longer time windows, while category 14 has shorter ones. In the case of shorter time windows, hospital deliveries of low-priority evacuees decrease, while shelters become more attractive to receive low-priority evacuees due to their proximity to the fire zone. Finally, the distinction between categories 14 and 32 is in the fire origin. In category 14, as well as all previously analyzed categories, the fire starts from the north. In category 32, however, the fire originates in the south, allowing more time to evacuate northern assembly areas. Figure 7 shows that in category 32, hospital 4 gains priority once again as it is located close to the northern part of the fire zone.

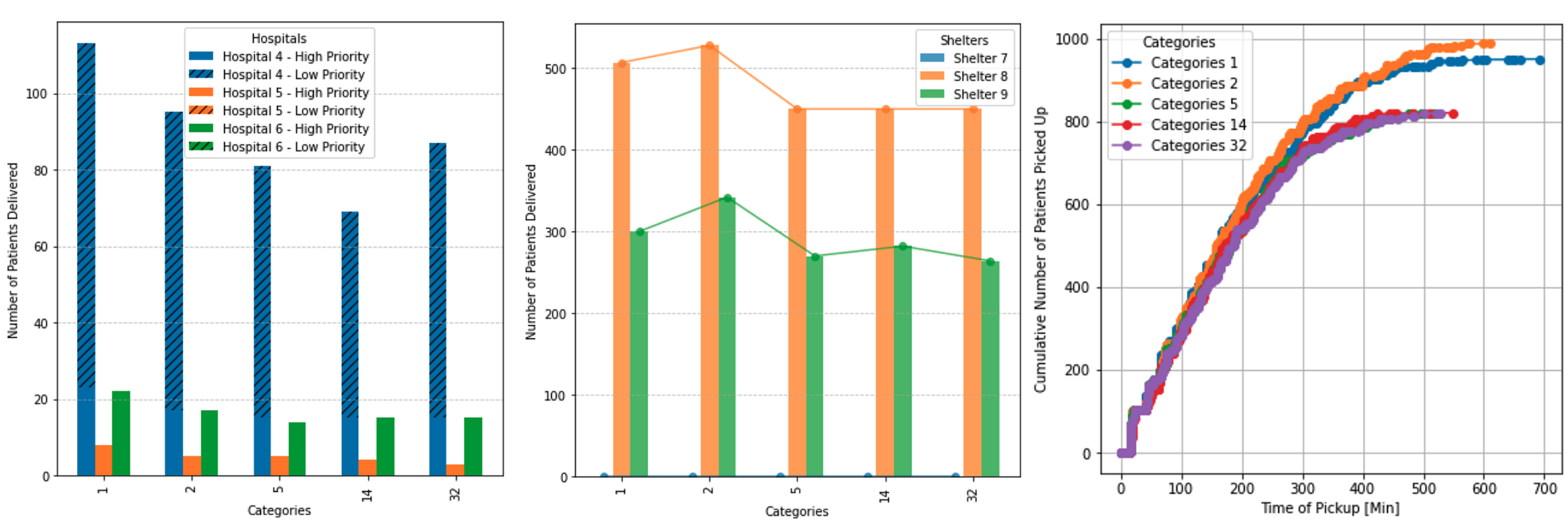

Figure 7: Results comparison in distinct categories for the case study

Table 7: Comparison of distinct categories

| Characteristic | | Category 1 | 2 | 5 | 14 | 32 |
|---|---|---|---|---|---|---|
| Fire Origin | | N | N | N | S | S |
| Total Population | | 950 (H) | 987 (H) | 821 (L) | 820 (L) | 819 (L) |
| P1 Population | | 53 (H) | 39 (L) | 34 (L) | 34 (L) | 33 (L) |
| Avg. Time Window | | 708 (Lo) | 684 (Lo) | 671 (Lo) | 473 (Sh) | 496 (Sh) |
| Hospital Drop-offs | H4 | 113 | 95 | 81 | 69 | 87 |
| | H5 | 8 | 5 | 5 | 4 | 3 |
| | H6 | 22 | 17 | 14 | 15 | 15 |
| Shelter Drop-offs | S7 | 0 | 0 | 0 | 0 | 0 |
| | S8 | 507 | 528 | 450 | 450 | 450 |
| | S9 | 300 | 342 | 270 | 282 | 264 |
| Evacuation time for given % of the population [min] | 50% | 158 | 157 | 142 | 143 | 144 |
| | 70% | 242 | 243 | 211 | 217 | 220 |
| | 100% | 646 | 574 | 482 | 509 | 500 |

## 7. Discussion

Timely evacuation is one of the most critical operational challenges in wildfire disaster management, particularly in countries such as Canada and the United States. Each year, vast areas of forest are devastated by wildfires, causing severe environmental damage, destroying wildlife habitats, and threatening human lives. During natural disasters like wildfires, most individuals can follow established protocols and evacuate to safe locations. However, vulnerable people, such as hospital patients, elderly individuals in long-term care facilities, and people with disabilities, cannot evacuate independently and require specialized transportation.

To address this issue, we formulate the supported-evacuation planning problem using a two-stage stochastic (TSS) programming approach. The first stage focuses on two critical resource management problems: (1) determining the required number of evacuation vehicles of each type and (2) selecting the number and locations of temporary shelters to receive low-priority evacuees. The second stage involves optimizing routing decisions, including pickups and drop-offs. To enhance model flexibility and realism, vehicles are allowed to perform split pickups and deliveries. Additionally, they do not have to return to their original dispatch location after each trip. Insights about minimum level of resources required are integrated into the solution methodology as logic-based inequalities. Specifically, when shelters offer similar services to different groups of patients, varying only in location and capacity, it is advantageous to activate as few shelters as possible to maximize capacity utilization. Since one of the objectives of activating shelters is to minimize transportation costs, priority should be given to shelters located at the minimum weighted average distance from assembly areas housing a larger number of low-priority evacuees.

We propose two TSS models based on different assumptions regarding the impact of wildfires and congestion on routes. In the absence of partial route disruptions, vehicles of the same type can be treated as independent and identical resources. This assumption allows for the development of an efficient solution methodology that leverages sequential decomposition techniques to simplify the original problem into more manageable sub-problems. Our results show that the proposed heuristics, used as benchmark policies, perform well, particularly in routing decisions. However, when considering partial route disruption or congestion, vehicles can no longer be treated as independent resources. In such cases, evacuation vehicles benefit from cooperative strategies. To account for that, we made minor modifications to our solution methodology to efficiently solve complex supported-evacuation instances within a reasonable time frame (ranging from a few minutes to less than an hour), significantly outperforming the direct solution provided by off-the-shelf solvers and the heuristic alternatives.

Building on this foundation, we developed a decomposition-based solution methodology to address computational complexities in large instances. This approach partitions a large, difficult-to-solve problem into multiple smaller, more tractable sub-problems while maintaining logical consistency. By relaxing linking constraints, we decomposed the original model by vehicle. Leveraging concepts from Logic-Based Benders decomposition, we divided the model into a restricted master problem (containing only first-stage variables) and $|S|$ scenario-based sub-problems. To accelerate convergence to the optimal solution, we introduced a neighborhood search through Benders combinatorial cuts. The resulting approach effectively reduces computational time and improves solution quality in terms of optimality gap and objective function.

Overall, the proposed framework offers a robust and efficient solution for wildfire evacuation planning for vulnerable populations. The methodology developed in this paper not only improves computational efficiency but also provides actionable strategies for real-world disaster response scenarios.

Despite our efforts to develop the most realistic model and address a problem that has received limited attention in the literature, the proposed approach has some limitations. One of our assumptions is that at the beginning of the evacuation, the number of patients in each assembly area is known and remains unchanged over time. However, this assumption might not be realistic as it is possible that people arrive at the assembly areas at different times. This dynamic arrival of people at assembly areas is an unexplored topic in supported-evacuation literature that could be pursued as a future research direction. While our original formulation of the problem does not allow changes in the values of parameter $Pa^{s}_{ih}$, our sequential decomposition approach allows for updates to this parameter which is the RHS of Constraint (43). However, the solution methodology would need to change to incorporate different scenarios regarding the number of new arrivals at the assembly areas to identify an optimum number of vehicles to be activated.

## Acknowledgments

This research was partially supported by the Natural Sciences and Engineering Research Council of Canada (NSERC) [Grants RGPIN-2018-05225 and RGPIN-2020-04301] and by the National Research Council of Canada (NRC) [Grant 18122020].

## References

Beloglazov, A., Almashor, M., Abebe, E., Richter, J., & Steer, K. C. B. (2016). Simulation of wildfire evacuation with dynamic factors and model composition. *Simulation Modelling Practice and Theory*, 60, 144–159.

Benders, J. (1962). Partitioning procedures for solving mixed-variables programming problems. *Numerische Mathematik*, 4(1), 238–252.

Carta, F., Zidda, C., Putzu, M., Loru, D., Anedda, M., & Giusto, D. (2023). Advancements in forest fire prevention: A comprehensive survey. *Sensors*, 23(14), 6635.

Clark, P. E., Porter, B. A., Pellant, M., Dyer, K., & Norton, T. P. (2023). Evaluating the efficacy of targeted cattle grazing for fuel break creation and maintenance. *Rangeland Ecology & Management*, 89, 69–86.

Codato, G. & Fischetti, M. (2006). Combinatorial Benders' cuts for mixed-integer linear programming. *Operations Research*, 54(4), 756–766.

Cova, T. J., Dennison, P. E., & Drews, F. A. (2011). Modeling evacuate versus shelter-in-place decisions in wildfires. *Sustainability*, 3(10), 1662–1687.

Dennison, P. E., Cova, T. J., & Mortiz, M. A. (2007). WUIVAC: A wildland-urban interface evacuation trigger model applied in strategic wildfire scenarios. *Natural Hazards*, 41, 181–199.

Dhall, A., Dhasade, A., Nalwade, A., VK, M. R., & Kulkarni, V. (2020). A survey on systematic approaches in managing forest fires. *Applied Geography*, 121, 102266.

Flores, I., Ortuño, M. T., & Tirado, G. (2023). A goal programming model for early evacuation of vulnerable people and relief distribution during a wildfire. *Safety Science*, 164, 106117.

Flores, I., Ortuño, M. T., Tirado, G., & Vitoriano, B. (2020). Supported evacuation for disaster relief through lexicographic goal programming. *Mathematics*, 8(4), 648.

Gan, H.-S., Richter, K.-F., Shi, M., & Winter, S. (2016). Integration of simulation and optimization for evacuation planning. *Simulation Modelling Practice and Theory*, 67, 59–73.

Geoffrion, A. M. (1972). Generalized Benders decomposition. *Journal of Optimization Theory and Applications*, 10, 237–260.

Goerigk, M., Deghdak, K., & Heßler, P. (2014). A comprehensive evacuation planning model and genetic solution algorithm. *Transportation Research Part E: Logistics and Transportation Review*, 71, 82–97.

Gwynne, S. M., Ronchi, E., Wahlqvist, J., Cuesta, A., Gonzalez Villa, J., Kuligowski, E. D., Kimball, A., Rein, G., Kinateder, M., Benichou, N., et al. (2023). Roxborough park community wildfire evacuation drill: Data collection and model benchmarking. *Fire Technology*, 59(2), 879–901.

Hooker, J. N. & Ottosson, G. (2003). Logic-based Benders decomposition. *Mathematical Programming*, 96, 33–60.

Kaiser, E. I., Hess, L., & Palomo, A. B. P. (2012). An emergency evacuation planning model for special needs populations using public transit systems. *Journal of Public Transportation*, 15(2), 45–69.

Kamyabniya, A. (2022). *Integrated and Coordinated Relief Logistics Planning Under Uncertainty for Relief Logistics Operations*. University of Ottawa.

Kamyabniya, A., Sauré, A., Salman, F. S., Bénichou, N., & Patrick, J. (2024). Optimization models for disaster response operations: A literature review. *OR Spectrum*, 46, 737–783.

Laporte, G. & Louveaux, F. V. (1993). The integer L-shaped method for stochastic integer programs with complete recourse. *Operations Research Letters*, 13(3), 133–142.

McCaffrey, S., Wilson, R., & Konar, A. (2018). Should I stay or should I go now? or should I wait and see? Influences on wildfire evacuation decisions. *Risk Analysis*, 38(7), 1390–1404.

McLennan, J., Ryan, B., Bearman, C., & Toh, K. (2019). Should we leave now? Behavioral factors in evacuation under wildfire threat. *Fire Technology*, 55, 487–516.

National Fire Protection Association (2024). *WUINITY: A platform for the simulation of wildland-urban interface fire evacuation*. https://www.nfpa.org/education-and-research/research/fire-protection-research-foundation/projects-and-reports/wuinity-a-platform-for-the-simulation-of-wildlandurban-interface-fire-evacuation. [Accessed: August 25, 2024].

National Forestry Database (2024). *National Fire Database*. https://cwfis.cfs.nrcan.gc.ca/ha/nfdb. Accessed: August 24, 2024

Nur, A. S., Kim, Y. J., Lee, J. H., & Lee, C.-W. (2023). Spatial prediction of wildfire susceptibility using hybrid machine learning models based on support vector regression in Sydney, Australia. *Remote Sensing*, 15(3), 760.

Rahmaniani, R., Crainic, T. G., Gendreau, M., & Rei, W. (2017). The Benders decomposition algorithm: A literature review. *European Journal of Operational Research*, 259(3), 801–817.

Shahparvari, S. & Abbasi, B. (2017). Robust stochastic vehicle routing and scheduling for bushfire emergency evacuation: An Australian case study. *Transportation Research Part A: Policy and Practice*, 104, 32–49.

Shahparvari, S., Abbasi, B., & Chhetri, P. (2017). Possibilistic scheduling routing for short-notice bushfire emergency evacuation under uncertainties: An Australian case study. *Omega*, 72, 96–117.

Stasiewicz, A. M. & Paveglio, T. B. (2021). Preparing for wildfire evacuation and alternatives: Exploring influences on residents' intended evacuation behaviors and mitigations. *International Journal of Disaster Risk Reduction*, 58, 102177.

Toledo, T., Marom, I., Grimberg, E., & Bekhor, S. (2018). Analysis of evacuation behavior in a wildfire event. *International Journal of Disaster Risk Reduction*, 31, 1366–1373.

Tymstra, C., Stocks, B. J., Cai, X., & Flannigan, M. D. (2020). Wildfire management in canada: Review, challenges and opportunities. *Progress in Disaster Science*, 5, 100045.

Weise, C. L., Brussee, B. E., Coates, P. S., Shinneman, D. J., Crist, M. R., Aldridge, C. L., Heinrichs, J. A., & Ricca, M. A. (2023). A retrospective assessment of fuel break effectiveness for containing rangeland wildfires in the sagebrush biome. *Journal of Environmental Management*, 341, 117903.

Zaidi, A. (2023). Predicting wildfires in algerian forests using machine learning models. *Heliyon*, 9(7), e18064.

Zhang, J., Zhang, M., & Li, G. (2021). Multi-stage composition of urban resilience and the influence of pre-disaster urban functionality on urban resilience. *Natural Hazards*, 107, 447–473.

Zhao, X., Xu, W., Ma, Y., & Hu, F. (2015). Scenario-based multi-objective optimum allocation model for earthquake emergency shelters using a modified particle swarm optimization algorithm: A case study in Chaoyang District, Beijing, China. *PLOS One*, 10(12), e0144455.

Zhao, X., Xu, Y., Lovreglio, R., Kuligowski, E., Nilsson, D., Cova, T. J., Wu, A., & Yan, X. (2022). Estimating wildfire evacuation decision and departure timing using large-scale GPS data. *Transportation Research Part D: Transport and Environment*, 107, 103277.

Zhou, S. & Erdogan, A. (2019). A spatial optimization model for resource allocation for wildfire suppression and resident evacuation. *Computers & Industrial Engineering*, 138, 106101.

## Appendix A. Proof of the Theorem

$$NV_k = \begin{cases} \min\left\{AV_k, \left\lceil \frac{\sum_{h\in H_k} MP_h}{CV_k \times \frac{Max_T}{2}} \right\rceil \right\} & \text{if } k \in \{1,2\} \\ \min\left\{AV_k, \left\lceil \frac{RP}{CV_k \times \frac{Max_T}{2}} \right\rceil \right\} & \text{if } k = \{3\} \text{ and } \sum_{k'\in K\setminus\{3\}} NV_{k'} > \sum_{k\in K\setminus\{3\}} AV_{k'} \\ 0 & \text{if } k = \{3\} \text{ and } \sum_{k'\in K\setminus\{3\}} NV_{k'} \le \sum_{k\in K\setminus\{3\}} AV_{k'} \end{cases} \tag{A.1}$$

where $Max_T$ is the maximum number of stops a ground vehicle can make, $\lceil\cdot\rceil$ denotes the ceiling function, which rounds a number up to the nearest integer, and $RP$ is the maximum total number of patients, across all scenarios, that ground vehicles are not expected to evacuate.

$$RP = \max\left\{0, \sum_{k\in K\setminus\{3\}} \left(\sum_{h\in H_k} MP_h - NV_k \times CV_k \times \frac{Max_T}{2}\right)\right\} \tag{A.2}$$

**Theorem 1.** *$NV_k$ is a valid lower bound on the number of vehicles of type $k \in \{1,2\}$ for full evacuation.*

**Proof** A valid lower bound on the number of vehicles of type $k$ is obtained when the number of pickups is set to its hypothetical upper bound. The maximum total number of patients a vehicle of type $k$ can pick up occurs when the vehicle makes $\frac{Max_T}{2}$ end-to-end trips between assembly areas and medical facilities without split pickup or delivery, and in each trip, the vehicle is loaded to its full capacity. Therefore, $CV_k \times \frac{Max_T}{2}$ is a valid upper bound on the number of patients picked up. Considering that the maximum total number of patients vehicles of type $k$ may need to evacuate is $\sum_{h\in H_k} MP_h$, complete evacuation of these people cannot be accomplished with fewer than $\left\lceil \frac{\sum_{h\in H_k} MP_h}{CV_k \times \frac{Max_T}{2}} \right\rceil$ vehicles. To account for the actual availability of vehicles of type $k$ in Equation (A.1), we take the minimum between this lower bound and $AV_k$. Thus, $NV_k$ is a valid lower bound on the number of vehicles of type $k$. □

## Appendix B. Algorithmic Description of the Solution Approach

**Algorithm 1** Pseudocode of the solution methodology

1: Set $r = 1$ and $R$ = maximum iterations
2: Solve RMP and define $S^0 = (\vec{X}^0, \vec{Y}^0)$
3: Solve SPs to obtain $\theta_v^s \quad \forall s, v$
4: Calculate $Z_0 = Z_{RMP}^0 + \sum_{s,v} \theta_v^s$
5: Best_sol = ($S^0$, second-stage variables)
6: Define Temporary Cut (TC) to make $S^0$ infeasible
7: $BS$ = [ ]"list of sub-optimal solutions"
8: **while** $r \leq R$ **OR** "stopping condition" is not met **do**
9: Add combinatorial cuts based on $S^{(r-1)}$
10: Add perturbation cuts
11: Add Temporary Cut (TC) to make $S^{(r-1)}$ infeasible
12: Add a cut to avoid visiting solutions in $BS$ (if the list is not empty)
13: Solve the RMP and SPs and define $S'^r = (\vec{X}^r, \vec{Y}^r)$
14: **if** RMP gets infeasible **then** "stopping condition" is met
15: Break
16: **end if**
17: Calculate $Z_r = Z_{RMP}^r + \sum_{s,v} \theta_v^s$
18: **if** $Z_r > Z_{r-1}$ **then** "*Bad Move*"
19: Add $S'^r$ to $BS$
20: $S^r = S^{(r-1)}$
21: Update UBs for the first-stage variable
22: **else**"*Good move*"
23: Define $S^r = S'^r$
24: Update Best_sol = ($S^r$, second-stage variables)
25: Update LBs for the first-stage variable
26: $r = r + 1$
27: **end if**
28: **end while**
29: Report Best_sol

## Appendix C. Practical Instance Results

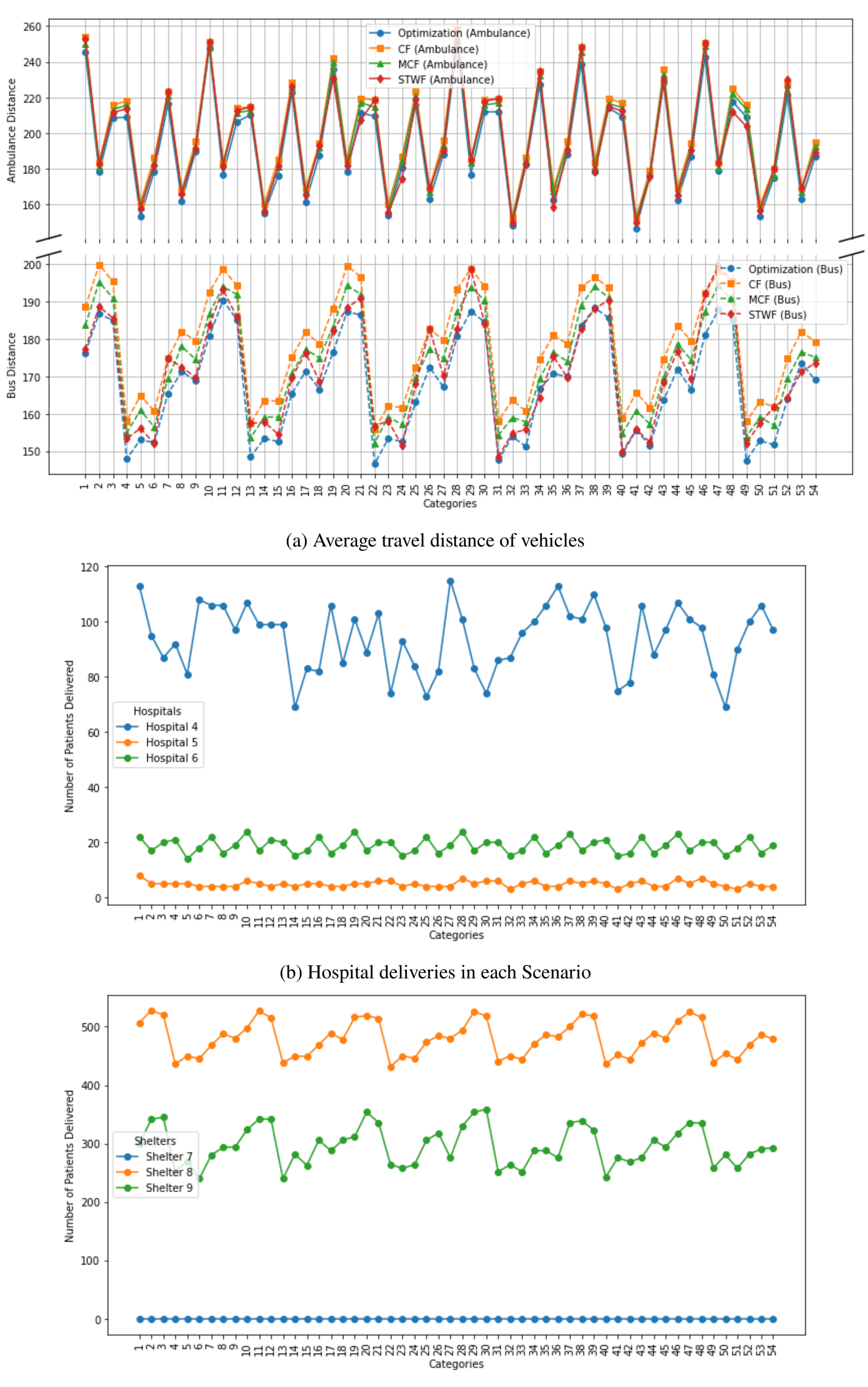

(a) Average travel distance of vehicles

(b) Hospital deliveries in each Scenario

(c) Shelter deliveries in each Scenario

Figure C.1: Overview of scenario results for the case study